\documentclass[11pt]{amsart}

\usepackage{euscript}
\usepackage{amsmath}
\usepackage{amsthm}
\usepackage{epsfig}
\usepackage{amssymb}

\usepackage{epic}

\DeclareMathOperator{\Hom}{{\rm Hom}}

\DeclareMathOperator{\coker}{{\rm coker}}

\DeclareMathOperator{\rank}{{\rm rank}}

\theoremstyle{plain}
\newtheorem{theorem}{Theorem}[subsection]

\theoremstyle{definition}

\newtheorem{para}[theorem]{}

\numberwithin{equation}{section}

\title{Lattice cohomology of normal surface singularities}

\author{Andr\'as N\'emethi}
\address{Alfr\'er R\'enyi Institute of Mathematics, 1053 Budapest, Re\'altanoda u. 13-15,
Hungary} \email{nemethi@renyi.hu}
\urladdr{http://www.renyi.hu/\textasciitilde nemethi}
\thanks{The author is partially supported by OTKA Grants.}



\date{}

\begin{document}

\maketitle

\begin{center}
 {\em Dedicated to Professor Heisuke Hironaka on his 77th birthday.}
\end{center}

\pagestyle{myheadings} \markboth{{\normalsize Andr\'as
N\'emethi}}{ {\normalsize Lattice cohomology}}

\newcommand{\pic}{{\rm Pic}^V_{top}(\Sigma)}
\newcommand{\picm}{{\rm Pic}_{top}(M)}
\newcommand{\pict}{{\rm Pic}_{top}(\Pi_{\vec{a}})}
\newcommand{\one}{\hat{\bf 1}}
\newcommand{\hD}{\hat{D}}
\newcommand{\eps}{\epsilon}
\newcommand{\ssw}{{\bf sw}}

\newcommand{\et}{\mathcal{T}}
\newcommand{\bS}{{\mathbb S}}
\newcommand{\bma}{\mbox{\boldmath$a$}}
\newcommand{\bmb}{\mbox{\boldmath$b$}}
\newcommand{\bmc}{\mbox{\boldmath$c$}}
\newcommand{\bme}{\mbox{\boldmath$e$}}
\newcommand{\bmi}{\mbox{\boldmath$i$}}
\newcommand{\bmj}{\mbox{\boldmath$j$}}
\newcommand{\bmv}{\mbox{\boldmath$v$}}
\newcommand{\bmk}{\mbox{\boldmath$k$}}
\newcommand{\bmm}{\mbox{\boldmath$m$}}
\newcommand{\bms}{\mbox{\boldmath$s$}}
\newcommand{\bSW}{\mbox{\boldmath$SW$}}
\newcommand{\bmf}{\mbox{\boldmath$f$}}
\newcommand{\bmg}{\mbox{\boldmath$g$}}
\newcommand{\gq}{{\mathfrak q}}
\newcommand{\xo}{o}

\newcommand{\gh}{g}

\newcommand{\lp}{{l}}
\newcommand{\ev}{\varepsilon}
\newcommand{\tx}{\tilde{X}}
\newcommand{\tz}{\tilde{Z}}
\newcommand{\call}{{\mathcal L}}
\newcommand{\calm}{{\mathcal M}}
\newcommand{\calx}{{\mathcal X}}
\newcommand{\calo}{{\mathcal O}}
\newcommand{\calt}{{\mathcal T}}
\newcommand{\cali}{{\mathcal I}}
\newcommand{\calj}{{\mathcal J}}
\newcommand{\calC}{{\mathcal C}}
\newcommand{\calS}{{\mathcal S}}
\newcommand{\calQ}{{\mathcal Q}}
\newcommand{\calF}{{\mathcal F}}

\newcommand{\cs}{\langle \chi_0\rangle}

\newcommand{\cc}{\bar{C}}
\newcommand{\vp}{\varphi}

\def\mmod{\mbox{mod}}
\let\d\partial
\def\EE{\mathcal E}
\newcommand{\cC}{\EuScript{C}}
\def\C{\mathbb C}
\def\Q{\mathbb Q}
\def\R{\mathbb R}
\def\bS{\mathbb S}
\def\bH{\mathbb H}
\def\bL{\mathbb L}
\def\Z{\mathbb Z}
\def\N{\mathbb N}\def\bn{\mathbb N}
\def\bp{\mathbb P}
\def\eop{$\hfill\square$}
\def\bif{(\, , \,)}
\def\coker{\mbox{coker}}
\def\im{{\rm Im}}

\newcommand{\Gammma}{{G}}

\newcommand{\p}{{tr}}
\newcommand{\Y}{A_{\alpha_1,\ldots,\alpha_p}}
\newcommand{\ap}{\tilde{A}^p}
\newcommand{\supno}{\stackrel{\supset}{\not=}}
\newcommand{\subno}{\stackrel{\subset}{\not=}}
\newcommand{\zbz}{Z_{B_0}}
\newcommand{\zbe}{Z_{B_1}}
\newcommand{\zbk}{Z_{B_2}}
\newcommand{\zm}{Z_{B_m}}
\newcommand{\zme}{Z_{B_{m-1}}}
\newcommand{\zmk}{Z_{B_{m-2}}}
\newcommand{\zi}{Z_{B_i}}
\newcommand{\zj}{Z_{B_j}}
\newcommand{\bz}{{\mathbb Z}}
\newcommand{\g}{E}
\newcommand{\zn}{Z_{min}}
\newcommand{\zmin}{Z_{min}}
\newcommand{\zk}{Z_K}
\newcommand{\co}{{\mathcal O}}
\newcommand{\x}{$(X,0)$\ }
\newcommand{\cl}{{\mathcal L}}
\newcommand{\ci}{{\mathcal I}_{D'}}
\newcommand{\y}{S. S.-T. Yau\ }
\newcommand{\no}{\noindent}
\newcommand{\bfc}{{\mathbb C}}
\newcommand{\bfq}{{\mathbb Q}}
\newcommand{\cale}{{\mathcal E}}
\newcommand{\calw}{{\mathcal W}}
\newcommand{\calv}{{\mathcal V}}
\newcommand{\calP}{{\mathcal P}}
\newcommand{\calI}{{\mathcal I}}
\newcommand{\cala}{{\mathcal A}}
\newcommand{\calr}{{\mathcal R}}
\newcommand{\calG}{{\mathcal G}}
\newcommand{\etan}{\eta(f;N)}
\newcommand{\etae}{\eta(f;1)}
\newcommand{\calp}{{\mathcal P}}
\newcommand{\f}{$f:({\bf C}^{n+1},0)\to ({\bf C},0)$\ }
\newcommand{\be}{B_{\epsilon}}
\newcommand{\se}{S_{\epsilon}}
\newcommand{\bc}{{\mathbb C}}
\newcommand{\bez}{B_{\epsilon_0}}
\newcommand{\br}{{\mathbb R}}
\newcommand{\bq}{{\mathbb Q}}
\newcommand{\sez}{S_{\epsilon_0}}
\newcommand{\ep}{\epsilon}
\newcommand{\vs}{\vspace{3mm}}
\newcommand{\fn}{f+z^N}
\newcommand{\si}{\sigma}
\newcommand{\Gammas}{S}

\begin{abstract}
For any negative definite plumbed 3--manifold $M$ we construct 
from its plumbed graph a graded $\Z[U]$--module. This, for rational
homology spheres, conjecturally equals the Heegaard--Floer homology of 
Ozsv\'ath and Szab\'o, but it has even more structure. If $M$ is a 
complex singularity link then the normalized Euler-characteristic can
be compared with the analytic invariants. The Seiberg--Witten
Invariant Conjecture  of \cite{[51],Line} is discussed in the light of
this new object. 
\end{abstract}


{\small

\section{Introduction}

The article is a symbiosis of singularity theory and low-dimensional
topology. Accordingly, it is preferable to separate its goals 
in two cathegories. 

\vspace{2mm}

From the point of view of {\em 3--dimensional topology}, the
article contains the following main result. For every negative
definite plumbed 3--manifold it constructs a graded
$\Z[U]$--module from the combinatorics of the plumbing graph.
This for rational homology spheres conjecturally equals the
Heegaard--Floer homology of Ozsv\'ath and Szab\'o. In fact, it has
more structure (e.g. instead of a $\Z_2$, odd/even grading it has a 
$\Z$ grading like a usual homology, see (\ref{HF7})(c)). 
The existence of these extra 
structures for arbitrary 3--manifolds might be an interesting subject
for further investigation. 

\vspace{2mm}

The motivations and aims from the pointof view of 
{\em singularity theory} are the following.  

In \cite{[51]} L. Nicolaescu and the author formulated a
conjecture which relates the geometric genus of a complex analytic
normal surface singularity $(X,0)$ --- whose link $M$ is a
rational homology sphere --- with the Seiberg-Witten invariant of
$M$ associated with the canonical $spin^c$--structure. The
conjecture generalized a conjecture of Neumann and Wahl \cite{NW}
which formulated the relationship for complete intersection
singularities with integral homology sphere links. The conjecture
\cite{[51]}  was verified in different cases, see
\cite{BN,NOSZ,[51],[52],[55],NO}.

Since the Seiberg-Witten theory provides a rational number for
{\em any} $spin^c$--structure, it was a natural challenge to
search for a complete set of conjecturally valid identities, which
involve all $spin^c$--structures. The preprint \cite{Line}
proposed such identities, connecting the sheaf--cohomology of
holomorphic line bundles associated with the analytic type of the
singularity with the Seiberg--Witten invariants of the link. The
identities were supported by a proof valid for rational
singularities.

But, a few months later,  \cite{SI} appeared with a list of
counterexamples. This posed a lot of questions: what kind of
guiding principles were wrongly interpreted in the original
conjectures? How can one `correct' them?

The present manuscript aims to answer some of them.

First, let us recall in short the original conjecture
(for canonical $spin^c$--structure). One fixes a topological type
(identified by a rational homology sphere link) and considers the
Seiberg--Witten invariant of this link (normalized with a certain
invariant $K^2+s$, see below). About this the conjecture predicted
two things: First, that it is an upper bound for the geometric
genus of all the possible analytic structures supported by the
fixed topological type. Second, that this bound is optimal, and it
is realized by all $\Q$--Gorenstein analytic structures.

Well, both expectations were wrong, but the nature of the two
errors are completely different. Regarding the second part, the
`Seiberg--Witten invariant identity', the error can be localized
easily. Indeed, the conjecture was over--optimistic: the identity
is not valid for {\em every} $\Q$--Gorenstein singularity.
Nevertheless, it is proved for large classes of singularities, and
we expect that the list will be continued. Hence, the form of the
identity shouldn't be modified, just we expect its validity for a
{\em subclass} of $\Q$--Gorenstein singularities. At this moment,
it is hopeless to identify exactly this subclass, nevertheless, in
\cite{BN} it has a description exclusively in terms of the
analytic structure
--- independently of the Seiberg--Witten theory; and \cite{Okuma2}
suggests that it can be identified by some vanishing
properties.

In fact, we were more concern about the inequality part:  Laufer
type computation sequences identified a possible topological upper
bound for the geometric genus, which in all cases explicitly
analysed (at the time of \cite{[51],Line}) coincided with the
Seiberg--Witten invariant, and the computation sequence
technique resonated perfectly with the theory of Ozsv\'ath and
Szab\'o from \cite{OSzP}. Then, which part of this line of
argument fails in general? The present article gives the following
answer: There exists a cohomology theory $\{\bH^q\}_{q\geq 0}$,
such that its normalized Euler--characteristic  (conjecturally)
equals the Seiberg--Witten invariant. On the other hand, its
$0^{th}$ normalized `Betti--number' (or invariants related with
it) serves as topological upper bound for the geometric genus (and
fits with computation sequence constructions). In simpler cases
(e.g. for rational, elliptic or star--shaped resolution graphs)
one has a vanishing $\bH^q=0$ for all $q\geq 1$, hence the
Seiberg--Witten invariant was able to serve as an upper bound.
But, in general, this is not the case: the geometric genus of
those analytic structures for which the `Seiberg--Witten invariant
identity' holds, are not extremal.

\vspace{2mm}

The article starts with the construction of this cohomology
theory: the {\em lattice cohomology}. Here, we do not restrict
ourselves to the rational homology sphere case. 
The construction provides from the plumbing graph of the link $M$
(or, from the associated intersection lattice) a graded
$\Z[U]$--module $\bH^q(M,\sigma)$ for each $q\geq 0$, and for all
torsion $spin^c$--structures $\sigma$ of $M$.

We emphasize and exemplify more the case $\bH^0$. $\bH^0$, as a
combinatorial $\Z[U]$--module associated with the link, is not new
in the literature: it was considered by Ozsv\'ath and Szab\'o in
\cite{OSzP} in Heegaard--Floer homology computations of some
special plumbed 3--manifolds (under the notation $\bH^+$). Later,
in \cite{NOSZ}, the author  computed $\bH^0$ for a larger class of
3--manifolds (`almost rational' graph--manifolds). In the present
article, in section 4, we prove similar characterization and
structure results for $\bH^*$ valid for rational, elliptic and
almost rational graphs. Moreover, we analyze examples with
$\bH^1\not=0$ too.
Section 5 connects $\bH^*$ with the Heegaard--Floer homology.

Section 6 deals with the theory of line bundles associated with
surface singularities. (It contains some parts from the
unpublished \cite{Line} and from the lecture notes \cite{trieste}.
Some similar $h^1$--computations for the case of rational
singularities were also found independently by T. Okuma
\cite{Okuma}.) In this section  we determine a topological upper
bound for the dimension of the sheaf--cohomologies of these line
bundles in terms of their Chern classes. The description sits in
$\bH^0$.

The last section 7 presents the `Seiberg--Witten invariant
conjecture' (the unmodified  conjectured identities), with
examples and more comments.

\vspace{2mm}

\noindent {\bf Acknowledgement.} The author was strongly
influenced by the articles of P\'eter Ozsv\'ath and Zolt\'an
Szab\'o, and benefited enormously from private conversations with
them. In fact, the idea of the existence of a cohomology theory
$\bH^*$ (for plumbed rational homology spheres) was formulated by
them, and  Conjecture (\ref{HF5}) is also part of their program.
The author also thanks the help of Andr\'as Stipsicz in different
aspects of Heegaard--Floer theory.

The first version of the manuscript was read by my colleagues
G\'abor Braun and Gyula Lakos, they provided valuable comments and
helped in finalization of some of the proofs.

\section{Preliminaries.}

\subsection{Negative definite plumbing graphs.}

\para\label{RES} Let $(X,0)$ be a complex analytic normal surface singularity
with link $M$. Fix a sufficiently small Stein representative $X$
of the germ $(X,0)$ and  let $\pi:\tilde{X}\to X$ be a {\em good}
resolution of the singular point $0\in X$. Let $E:=\pi^{-1}(0)$ be
the exceptional divisor with irreducible components $\{E_j\}_{j\in
{\mathcal J}}$ and write $\Gamma(\pi)$ for the dual resolution
graph associated with $\pi$.  Recall that $\Gamma(\pi)$ is
connected and the intersection matrix $I:=\{(E_j,E_i)\}_{j,i}$ is
negative definite. We write $e_j$ for $E_j^2$,  $g_j$ for the
genus of $E_j$ ($j\in {\mathcal J}$), and $g:=\sum_jg_j$.
Moreover, let $c$ be the number of independent cycles in (the
topological realization of) $\Gamma$. E.g., $c=0$ if and only if
\, $\Gamma(\pi)$ is a tree. The rank of $H_1(M,\Z)$ is
$c+2g$. Hence, $M$ is a rational homology sphere  (i.e.\
$H_1(M,\Q)=0$)  if and only if $g=c=0$.

\para Since  $\pi $ identifies $\partial
\tilde{X}$ with $M$, the  graph $\Gamma(\pi)$ can be viewed as a
plumbing graph and $M$ as the associated $S^1$-plumbed manifold.
In the sequel $\Gamma$ will denote either a good resolution graph
as above, or a negative definite plumbing graph of $M$. Similarly,
$\tilde{X} $ denotes either the space of a good resolution, or the
oriented 4-manifold obtained by plumbing disc-bundles
corresponding to $\Gamma$.

\subsection{The combinatorics of the plumbing.}

\para\label{lattice}{\bf Definition. The lattices $L$ and $L'$.}
The image of the boundary operator
$\partial: H_2(\tilde{X},M,\Z)\to H_1(M,\Z)$ is the
torsion subgroup $H$  of $H_1(M,\Z)$. The
exact sequence of $\Z$-modules
\begin{equation*}
0\to L \stackrel{i}{\to}  L'\to H\to 0 \tag{1}
\end{equation*}
will stand for the homological exact sequence
\begin{equation*}
0\to H_2(\tx,\Z)\to
H_2(\tx,M,\Z)\stackrel{\partial}{\longrightarrow} Tors(
H_1(M,\Z))\to 0,
\end{equation*}
(or for its Poincar\'e dual).  Here $L$ is freely generated by the
homology classes $\{E_j\}_{j\in \calj}$ and is equipped with the
intersection form $(\cdot,\cdot)$. For each $j$, consider a small
transversal disc $D_j$ in $\tx$ with $\partial D_j\subset
\partial \tx$. Then $L'$ is freely generated by the (relative
homology) classes $\{D_j\}_{j\in \calj}$.
Notice that  the morphism $i:L\to L'$ can be
identified with $L\to \Hom(L,\Z)$ given by $l\mapsto (l,\cdot)$.
The intersection form has a natural extension to $L_\Q=L\otimes
\Q$, and we will regard $\Hom(L,\Z)$ as a sub-lattice  of
$L_\Q$: $\alpha\in \Hom (L,\Z)$ corresponds with the unique
$l_\alpha\in L_\Q$ which satisfies $\alpha(l)=(l_\alpha,l)$ for
any $l\in L$.
Hence, the exact sequence (1) can be recovered completely from the
lattice $L$.

\para\label{2.3}{\bf Characteristic elements. $Spin^c$-structures.} The
set of characteristic elements are defined by
$$Char=Char(L):=\{k\in L': \, (k,x)+(x,x)\in 2\Z \ \mbox{for any $x\in L$}\}.$$
The unique rational cycle $K\in L'$ which satisfies the system (of
adjunction relations) $(K,E_j)=-(E_j,E_ j)- 2+2g_j$ for all $j$
is called the  {\em canonical cycle}. Then $Char=K+2L'$. There is
a natural action of $L$ on $Char$ by $l*k:=k+2l$ whose orbits are
of type $k+2L$. Obviously, $H$ acts freely and transitively on the
set of orbits by $[l']*(k+2L):=k+2l'+2L$.

If $\tilde{X}$ is a 4--manifold as above, then $H^2(\tilde{X},\Z)$
has no 2--torsion. Therefore, the first Chern class (of the
associated determinant line bundle) realizes an identification
between the $spin^c$--structures $Spin^c(\tilde{X})$ on
$\tilde{X}$ and $Char\subset L'= H^2(\tilde{X},\Z)$ (see e.g.
\cite{GS}, 2.4.16). On the other hand, the $spin^c$--structures on
$M$ form an $H_1(M,\Z)$ torsor. In the image of the restriction
$Spin^c(\tilde{X})\to Spin^c(M)$ are exactly those
$spin^c$--structures of $M$ whose Chern classes are restrictions
$L'\to H^2(M,\Z)=H_1(M,\Z)$, i.e. are torsion elements  sitting in
$H$. We call them  \emph{torsion}
structures, and we denote them by $Spin^c_t(M)$. One has
an identification of $Spin^c_t(M)$  with the set of $L$--orbits of
$Char$, and this identification is compatible with the
action of $H$ on both sets. In the sequel, we think about
$Spin^c_t(M)$ by this identification:  any torsion
$spin^c$-structure of $M$ will be represented by an orbit
$[k]:=k+2L\subset Char$. The \emph{canonical } $spin^c$--structure
(is torsion and) corresponds to $[K]$.

We write $\hat{H}$ for the Pontrjagin dual
$\Hom(H,S^1)$ of $H$. One has a natural isomorphism
$$\theta:H\to \hat{H}, \ \ \mbox{induced by} \ \   [l']\mapsto e^{2\pi i (l',\cdot )}.$$

\para\label{2.2p}{\bf Positive cones.} One can consider two types of
`positivity conditions'
 for rational cycles.  The first one is considered in $L$.
A cycle $x=\sum_j  r_jE_j\in  L_\Q$ is  called {\em effective},
denoted by $x\geq 0$, if $r_j\geq 0$  for all $j$.  Their
collection  is denoted by $L_{\Q,e}$, while $L_e' :=L_{\Q,e}\cap
L'$ and $L_e :=L_{\Q,e}\cap L$.

The second is the {\em numerical effectiveness} of the rational
cycles,
 i.e.\ positivity considered in $L'$.
We define $L_{\Q,ne}:=\{x\in L_\Q : (x,E_j)\geq  0 \ \mbox{for all
$j$}\}$. In fact, $L_{\Q,ne}$ is the positive cone in $L_\Q$
generated by $\{D_j\}_j$, i.e.\ it is exactly $\{\sum_j r_j D_j,
r_j\geq 0 \ \mbox{for all $j$}\}$. Since $I$ is negative definite,
all the entries of $D_j$ are {\em strictly}  negative. In
particular, $-L_{\Q,ne} \subset L_{\Q,e}$. Similarly as above,
write $L_{ne}:=L\cap L_{\Q,ne}$.

%

\para\label{2.4}{\bf Liftings.} We will
consider some `liftings' (set theoretical sections) of the element
of $H$ into $L'$. They correspond to the positive cones in $L_\Q$
considered in (\ref{2.2p}).

More precisely, for any $l'+L=h\in H$, let $l'_e(h)\in L'$ be the
unique minimal effective rational cycle in $L_{\Q,e}$ whose class
is $h$. Clearly, the set $\{l'_e(h)\}_{h\in H}$ is exactly
$Q:=\{\sum_j r_jE_j \in L'\,;\, 0\leq r_j<1\}$. 

Similarly, for any $h=l'+L$, the intersection $(l'+L)\cap
L_{\Q,ne}$ has a unique maximal element $l'_{ne}(h)$, and the
intersection $(l'+L)\cap (-L_{\Q,ne})$) has a unique minimal
element $\bar{l}'_{ne}(h)$
 (cf.\ \cite{NOSZ}, 5.4).  By their definitions
$\bar{l}'_{ne}(h)=-l'_{ne}(-h)$.


For some $h$, $\bar{l}'_{ne}(h)$ might be situated in $Q$, but, in
general, this is not the case. In general, the characterization of
all the elements $\bar{l}'_{ne}(h)$ is  not simple (see e.g.
\cite{NOSZ}).

\para\label{chifunction}{\bf The $\chi$-functions (Riemann-Roch formula).} For any
characteristic element $k\in Char$ one defines $$\chi_k:L'\to\Q\ \
\mbox{by} \ \ \chi_k(l'):=-(l',l'+k)/2.$$ Clearly,
$\chi_k(L)\subset \Z$. For the interpretation of $\chi_k$ as
(twisted) Riemann-Roch formula, consider the following. Let
$\tilde{X}$ be a resolution as in (\ref{RES}), and fix a
holomorphic line bundle $\call\in Pic(\tx)$, and write
$c_1(\call)=l'\in L'$ for its Chern class. 
Set $k:=K-2l'\in Char$. For any $l\in L$ with $l>0$ one defines
the sheaf $\calo_l:=\calo_{\tx}/\calo_{\tx}(-l)$ supported by $E$
(see e.g. \ref{CohCom}).
Consider the sheaf $\call\otimes \calo_l$  and let
$\chi(\call\otimes \calo_l)=h^0(\call\otimes
\calo_l)-h^1(\call\otimes \calo_l)$ be its (holomorphic)
Euler-characteristic. The Riemann-Roch theorem states that this
can be computed  combinatorially, namely
$$\chi(\call\otimes \calo_l)=\chi_k(l).$$

\section{The lattice cohomology associated with $L$.}

\subsection{Lattice cohomology associated with $\Z^s$ and a system
of weights.}

\para\label{a} We consider a free $\Z$--module, with a fixed basis
$\{E_j\}_j$, denoted by $\Z^s$. It is also convenient to fix
a total ordering of the index set ${\mathcal J}$, which in the
sequel will be denoted by $\{1,\ldots,s\}$.

Our goal is to define a graded $\Z[U]$--module associated with the
pair $(\Z^s, \{E_j\}_j)$ and a system of weights, which will be
introduces in (\ref{weight}). First we set some notations
regarding $\Z[U]$--modules.

\para\label{zu1} {\bf $\Z[U]$-modules.}  Consider the graded $\Z[U]$--module
$\et:=\Z[U,U^{-1}]$, and  (following \cite{OSzP}) denote by
 $\calt_0^+$ its quotient by the submodule  $U\cdot \Z[U]$.
This has a grading in such a way that $\deg(U^{-d})=2d$ ($d\geq
0$). Similarly, for any $n\geq 1$,  the quotient of $\Z\langle
U^{-(n-1)}, U^{-(n-2)},\ldots, 1,U,\ldots \rangle$ by $U\cdot
\Z[U]$ (with the same grading) defines the  graded module
$\calt_0(n)$. Hence, $\calt_0(n)$, as a $\Z$--module, is freely
generated by $1,U^{-1},\ldots,U^{-(n-1)}$, and has finite
$\Z$-rank $n$.

More generally, for any graded $\Z[U]$--module $P$ with
$d$--homogeneous elements $P_d$, and  for any  $r\in\Q$,   we
denote by $P[r]$ the same module graded (by $\Q$) in such a way
that $P[r]_{d+r}=P_{d}$. Then set $\calt^+_r:=\calt^+_0[r]$ and
$\calt_r(n):=\calt_0(n)[r]$. (Hence, for $m\in \Z$,
$\et_{2m}^+=\Z\langle U^{-m}, U^{-m-1},\ldots\rangle$.)

\para\label{complex} {\bf The cochain complex.}
$\Z^s\otimes \R$ has a natural cellular decomposition into cubes. The
set of zero--dimensional cubes is provided  by the lattice points
$\Z^s$. Any $l\in \Z^s$ and subset $I\subset {\mathcal J}$ of
cardinality $q$  defines a $q$--dimensional cube, which has its
vertices in the lattice points $(l+\sum_{j\in I'}E_j)_{I'}$, where
$I'$ runs over all subsets of $I$. On each such cube we fix an
orientation. This can be determined, e.g.,  by the order
$(E_{j_1},\ldots, E_{j_q})$, where $j_1<\cdots < j_q$, of the
involved base elements $\{E_j\}_{j\in I}$. The set of oriented
$q$--dimensional cubes defined in this way is denoted by $\calQ_q$
($0\leq q\leq s$).

Let $\calC_q$ be the free $\Z$--module generated by oriented cubes
$\square_q\in\calQ_q$. Clearly, for each $\square_q\in \calQ_q$,
the oriented boundary $\partial \square_q$ has the form
$\sum_k\varepsilon_k \, \square_{q-1}^k$ for some
$\varepsilon_k\in \{-1,+1\}$. Here, in this sum, we write only
those $(q-1)$--cubes which appear with non--zero coefficient.
These are called {\em faces} of $\square_q$.

It is clear that $\partial\circ\partial=0$. But, obviously, the
homology of the chain complex $(\calC_*,\partial)$
(or, of the cochain complex $(\Hom_{\Z}(\calC_*,\Z),\delta)$)
is not very
interesting: it is just the (co)homology of $\R^s$.
%
A more interesting (co)homology can be constructed as follows. For
this, we consider a set of compatible {\em weight
functions} $w_q:\calQ_q\to \Z$ ($0\leq q\leq s$).

\para\label{weight} {\bf Definition.} A set of functions
$w_q:\calQ_q\to \Z$  ($0\leq q\leq s$) is called a {\em set of
compatible weight functions}  if the following hold:

(a) For any integer $k\in\Z$, the set $w_0^{-1}(\,(-\infty,k]\,)$
is finite;

(b) for any $\square_q\in \calQ_q$ and for any of its faces
$\square_{q-1}\in \calQ_{q-1}$ one has $w_q(\square_q)\geq
w_{q-1}(\square_{q-1})$.

\noindent (In the sequel sometimes we will omit the index $q$ of
$w_q$.)

\vspace{2mm}

Assume that we already fixed a set of compatible weight functions
$\{w_q\}_q$. Then we set $\calF^q:=\Hom_{\Z}(\calC_q,\et^+_0)$.
Notice that $\calF^q$ is, in fact, a $\Z[U]$--module by
$(p*\phi)(\square_q):=p(\phi(\square_q))$ ($p\in \Z[U]$).
Moreover, $\calF^q$ has a $\Z$--grading: $\phi\in \calF^q$ is
homogeneous of degree $d\in\Z$ if for each $\square_q\in\calQ_q$
with $\phi(\square_q)\not=0$, $\phi(\square_q)$ is a homogeneous
element of $\et^+_0$ of degree $d-2\cdot w(\square_q)$. (In fact,
the grading is $2\Z$--valued; hence, the reader interested only in
the present construction may divide all the degrees by two.
Nevertheless, we prefer to keep the present form in our
presentation because of its resonance with the Heegaard--Floer
homology of the link.)

Next, we define $\delta_w:\calF^q\to \calF^{q+1}$. For
this, fix $\phi\in \calF^q$ and we show how $\delta_w\phi$ acts on
a cube $\square_{q+1}\in \calQ_{q+1}$. First write
$\partial\square_{q+1}=\sum_k\varepsilon_k \square ^k_q$, then
set
$$(\delta_w\phi)(\square_{q+1}):=\sum_k\,\varepsilon_k\,
U^{w(\square_{q+1})-w(\square^k_q)}\, \phi(\square^k_q).$$

\para\label{Lemma1} {\bf Lemma.} {\em $\delta_w\circ\delta_w=0$, i.e.
$(\calF^*,\delta_w)$ is a cochain complex.}

\begin{proof} With the obvious notations,
$(\delta_w^2\phi)(\square^j_{q+2})$ equals $$\sum_k\varepsilon
^j_kU^{w(\square_{q+2}^j)-w(\square^k_{q+1})}\,
\sum_l\varepsilon_l^k
U^{w(\square^k_{q+1})-w(\square^l_q)}\phi(\square^l_q)=
\sum_lU^{w(\square^j_{q+2})-w(\square^l_q)}\,\big(\,
\sum_k\varepsilon^j_k\varepsilon^k_l\,\big)\, \phi(\square^l_q).$$
But, for any $l$,  $\sum_k\varepsilon^j_k\varepsilon^k_l=0$ since
$\partial^2=0$.
\end{proof}

\para In fact, $(\calF^*,\delta_w)$ has a natural {\bf
augmentation} too. Indeed, set $m_w:=\min_{l\in \Z^s}w_0(l)$ and
choose  $l_w\in \Z^s$ such that $w_0(l_w)=m_w$. Then one defines
the $\Z[U]$--linear map
$$\epsilon_w:\et^+_{2m_w}\longrightarrow \calF^0$$ such that
$\epsilon_w (U^{-m_w-s})(l)$ is the class of $U^{-m_w+w_0(l)-s}$
in $\et^+_0$ for any integer $s\geq 0$.

\para\label{lemma2}{\bf Lemma.}  {\em $\epsilon_w$ is injective, and
$\delta_w\circ\epsilon_w=0$.}
\begin{proof} Since $\epsilon_w(U^{-m_w-s})(l_w)=U^{-s}$,
the injectivity is clear. Take $\square\in\calQ_1$ with
$\partial\square=a-b$. Then
$$(\delta_w\epsilon_w)(t)(\square)=U^{w(\square)-w(a)}\epsilon_w(t)(a)
-U^{w(\square)-w(b)}\epsilon_w(t)(b)=U^{w(\square)}t-U^{w(\square)}t=0.$$

\vspace*{-7mm}

\end{proof}

We invite the reader to verify that $\epsilon_w$ and $\delta_w$
are morphisms of $\Z[U]$--modules, and are homogeneous of degree
zero.

\para\label{def12}{\bf Definitions.} The homology of the cochain complex
$(\calF^*,\delta_w)$ is called the {\em lattice cohomology} of the
pair $(\R^s,w)$, and it is denoted by $\bH^*(\R^s,w)$. The
homology of the augmented cochain complex
$$0\longrightarrow\et^+_{2m_w}\stackrel{\epsilon_w}{\longrightarrow}
\calF^0\stackrel{\delta_w}{\longrightarrow}\calF^1
\stackrel{\delta_w}{\longrightarrow}\ldots$$ is called the {\em
reduced lattice cohomology} of the pair $(\R^s,w)$, and it is
denoted by $\bH_{red}^*(\R^s,w)$. If the pair $(\R^s,w)$ is clear
from the context, we omit it from the notation. Clearly, for any
$q\geq 0$, both $\bH^q$ and $\bH_{red}^q$ admit an induced graded
$\Z[U]$--module structure and $\bH^q=\bH^q_{red}$ for $q>0$.
Moreover, the $\Z$--grading of $\calF^q$ induces a $\Z$--grading
on $\bH^q$ and $\bH^q_{red}$; the homogeneous part of degree $d$
is denoted by $\bH^q_d$, or $\bH^q_{red,d}$.

It is easy to see that $\bH^*(\R^s,w)$ depends essentially on the
choice of $w$.

\para\label{lemma3}{\bf Lemma.} {\em One has a graded
$\Z[U]$--module isomorphism
$\bH^0=\et^+_{2m_w}\oplus\bH^0_{red}$.}

\begin{proof}
Consider the isomorphism $U^{-m_w}:\et^+_0\to \et^+_{2m_w}$. Then
define $r_w:\bH^0\to \et^+_{2m_w}$ by $r_w(\phi):=
U^{-m_k}\phi(l_w)$. Since $r_w\circ \epsilon_w =1$, the exact
sequence $0\to \et^+_{2m_w}\stackrel{\epsilon_w}{\longrightarrow}
\bH^0\to \bH^0_{red}\to 0$ \, splits.
\end{proof}

\para\label{rem} Next, we present another realization of the
modules $\bH^*$.

\para\label{Sn} {\bf Definitions.} For each $n\in \Z$,
define $S_n=S_n(w)\subset \R^s$ as the union of all
the cubes $\square_q$ (of any dimension) with $w(\square_q)\leq
n$. Clearly, $S_n=\emptyset$, whenever $n<m_w$. For any  $q\geq 0$, set
$$\bS^q(\R^s,w):=\oplus_{n\geq m_w}\, H^q(S_n,\Z).$$
Then $\bS^q$ is $\Z$ (in fact, $2\Z$)--graded, the
$d=2n$--homogeneous elements $\bS^q_d$ consist of  $H^q(S_n,\Z)$.
Also, $\bS^q$ is a $\Z[U]$--module; the $U$--action is given by
the restriction map $r_{n+1}:H^q(S_{n+1},\Z)\to H^q(S_n,\Z)$.
Namely,  $U*(\alpha_n)_n=(r_{n+1}\alpha_{n+1})_n$. Moreover, for
$q=0$, the fixed base--point $l_w\in S_n$ provides an augmentation
(splitting)
 $H^0(S_n,\Z)=
\Z\oplus \tilde{H}^0(S_n,\Z)$, hence an augmentation of the graded
$\Z[U]$--modules
$$\bS^0=\et^+_{2m_w}\oplus \bS^0_{red}=(\oplus_{n\geq m_w}\Z)\oplus (
\oplus_{n\geq m_w}\tilde{H}^0(S_n,\Z)).$$

\para\label{STR1} {\bf Theorem.} {\em

(a) There exists a graded $\Z[U]$--module isomorphism, compatible with the
augmentations:
$$\bH^*(\R^s,w)=\bS^*(\R^s,w).$$

(b) For any degree $d$, there exists an integer $N(d)\geq 0$, such that
$$\bH^*_{red,d}\cap im(U^{N(d)})=0.$$

(c) For any $n$ one has $U^{n-m_w+1}\bH^*_{2n}=0$. If there exists $N$
such that $S_n$ is
contractible for any $n\geq N$,  then  $U^{N-m_w}\bH^*_{red}=0$.}

\begin{proof} (a) Let $\calF_d^q$ be the set of $d=2n$--homogeneous elements
$\phi\in \calF^q$. Since $\delta_w$ is homogeneous of degree zero,
$(\calF^*_d,\delta_w)$ is a complex. Let $(\calC^*(S_n),\delta)$ be the usual
cochain
complex of $S_n$. Then the two complexes can be naturally identified.
Indeed, take $\phi\in \calF^q_d$. Then, for any $\square_q$, $\phi(\square_q)$
has the form $a_\phi(\square_q)U^{w(\square_q)-n}$.
Hence $a_{\phi}(\square_q)\in\Z$ is well--defined for any $q$--cube
$\square_q$ of $S_n$, and  the correspondence
$\phi\mapsto a_\phi$ realizes  the bijection $\calF^*_d\to \calC^*(S_n)$.

Since $\tilde{H_q}(\R^s,\Z)=0$, for any $n$ there exists $N$ such that
$\tilde{H}_q(S_n)\to \tilde{H}_q(S_{n+N})$ is trivial. (b) is the dual
statement of this. (c) follows from (a).

\end{proof}

\para\label{Fin} {\bf Remark.} Although 
$\bH^*_{red}(\R^s,w)$ has finite $\Z$--rank in any fixed
homogeneous degree, in general, it is not
finitely generated over $\Z[U]$.  E.g., set  $s=1$, and define $w_0$ by
$$w_0(-n)=w_0(n)=[n/2]+4\{n/2\} \ \ \mbox{for any $n\in\Z_{\geq 0}$},$$
where $[\,]$ and   $\{\,\}$ are the integral, respectively  the
fractional parts; and let $w_1$ on the segment $[n,n+1]$ take the
value $\max\{w_0(n),w_0(n+1)\}$. Then
$\bH^0_{red}=\oplus_{k\geq 1}\et_k(1)^2$.

\para\label{SSP} {\bf Restrictions.} Assume that $T\subset \R^s$ is a subspace
of $\R^s$ consisting of union of some cubes (from $\calQ_*$). Let
$\calC_q(T)$ be the free $\Z$--module generated by $q$--cubes of
$T$, $\calF^q(T)= \Hom_\Z(\calC_q(T),\et^+_0)$. Then
$(\calF^*(T),\delta_w)$ is a complex, whose homology will be
denoted by $\bH^*(T,w)$. It has a natural graded $\Z[U]$--module
structure.  The restriction map induces a natural graded
$\Z[U]$--module homogeneous homomorphism (of degree zero)
$$r^*:\bH^*(\R^s,w)\to \bH^*(T,w).$$

\subsection{Lattice cohomology associated with $\Gamma$ and $k\in Char$.}

\para\label{dEF1} We consider a graph $\Gamma$ as in \S 2
and we fix a characteristic element $k\in Char$. Notice that
$\Gamma$ automatically provides a free $\Z$--module $L=\Z^s$ with
a fixed bases $\{E_j\}_j$. Using $\Gamma$ and $k$, we define a set
of compatible weight functions $\{w_q\}_q$.

The definition reflects our effort to connect the topology of a
singularity--link (e.g. the lattice cohomology) with analytic
invariants. For more detailed motivation, see (\ref{El5}) and
(\ref{PCUB}).

For any $\gh\geq 0$ and  $n\geq 0$, let
$M^{\gh}(n)$ be the maximum of all possible dimensions of
sheaf--cohomologies $H^1(C,\call)$, where $C$ runs over all Riemann surfaces
of genus $\gh$
and $\call$ is a holomorphic line bundle on $C$ with holomorphic
Euler--characteristic $\chi(\call)=n$. (This number exists, in fact
$M^g(n)\leq g$.)

Now, we define $\{w_q\}_q$ as follows.  For $q=0$ we set
$w_0:=\chi_k$ (cf. \ref{chifunction}). Since the intersection form
is negative definite, (\ref{weight})(a) is satisfied.

Next, we define $w_1$. Consider a segment $S\in \calQ_1$ with
vertices $l$ and $l+E_j$ for some $l\in L$ and $j\in{\mathcal J}$.
We set
$$w_1(S):=\max\,\{\,\chi_k(l)\,,\,\chi_k(l+E_j)\,\}
+ M^{g_j}(\, |\chi_k(l)-\chi_k(l+E_j)|\,).$$
Finally, for any $\square_q\in \calQ_q$ ($q\geq 2$) set
$$w_q(\square_q):=\max\,\{\, w_1(S)\, :\, S \ \ \mbox{is a segment
of $\square_q$}\,\}.$$

\para\label{RELL} {\bf Examples.} (a) Assume that $g_j=0$ for all $j$.
Since  $M^0(n)=0$ for any $n\geq 0$, for any $q$
$$ w(\square_q)=\max\{\chi_k(v)\,:\, v\ \mbox{is a vertex of
$\square_q$}\}.$$

(b) Assume that $g_j\leq 1$ for any $j$. Since $M^1(n)=0$ for $n\geq 1$ and
$M^1(0)=1$, the definition of $w_1$ might be modified into
$$w_1(S)=\max\,\{\,\chi_k(l)\,,\,\chi_k(l+E_j)\, \, , \, \min\{
\,\chi_k(l)\,,\,\chi_k(l+E_j)\, \}\, +g_j\, \}.$$

(c) By a vanishing theorem, in general, $M^{\gh}(n)=0$ whenever $n\geq \gh$.

\para\label{DEF} {\bf Definition.} The  $\Z[U]$--modules $\bH^*(\R^s,w)$ and
$\bH^*_{red}(\R^s,w)$ obtained by these weight functions are
called the {\em lattice cohomologies} associated with the pair
$(\Gamma,k)$ and are denoted by $\bH^*(\Gamma,k)$, respectively
$\bH^*_{red}(\Gamma,k)$. We also write $m_k:=m_w= \min_{l\in L}
\chi_k(l)$.

\para\label{STR2} {\bf Theorem.} {\em
$\bH^*_{red}(\Gamma,k)$ is finitely generated over $\Z$.}

\begin{proof}
We start the proof with the following statement:

\vspace{1mm}

\noindent {\bf Fact:} {\em There exist $X\in L_e$  and an
increasing infinite sequence of cycles $\{x_i\}_{i\geq 0}$ with
$x_0=X$, such that

(a) $x_{i+1}=x_i+E_{j(i)}$ for some $j(i)$, $i\geq 0$,

(b) if $x_i=\sum_j m_{i,j}E_j$, then for all $j$,
$m_{i,j}$ tends to infinity as $i$ tends to infinity,

(c) $\chi_k(x_{i+1})-\chi_k(x_i)\geq g_{j(i)}$. }

\vspace{1mm}

{\em Similarly, there exists $Y\in L_e$  and an increasing
infinite sequence of cycles $\{y_i\}_{i\geq 0}$, with $y_0=Y$ and
similar properties as in (a)--(b), and (c)
$\chi_k(-y_{i+1})-\chi_k(-y_i)\geq g_{j(i)}$. }

\vspace{1mm}

Indeed, take a cycle $Z\in L$ such that $(Z,E_j)<0$ for any $j$.
Let $\{z_i\}_{i=0}^t$ be an increasing sequence with $z_0=0$,
$z_t=Z$, $z_{i+1}= z_i+E_{j(i)}$ ($0\leq i<t$). Then for $m$
sufficiently large, $X=mZ$, and the sequence $\{m'Z+z_i\}$ (where
$m'\geq m$ and $0\leq i< t$) works. A similar statement is valid for
$Y=mZ$ (and similar type of sequence) with  $m\gg 0$.

Fix $X,Y\in L_e$, such that $-Y\leq l_w\leq X$. Let
$T(-Y,X)=\{r\in \R^s\,:\, -Y\leq r\leq X\}$. $T(-Y,X)$ has a
natural cube--decomposition compatible with the decomposition of
$\R^s$, hence by (\ref{SSP}), one has a map $r^*_{-Y,X}:
\bH^*_{red}(\R^s,w)\to \bH^*_{red}(T(-Y,X),w)$.

Set  $X$ and $Y$ as in Fact; clearly we may assume that $-Y\leq
l_w\leq X$. We claim that  $r^*_{-Y,X}$ is an isomorphism. Indeed,
consider the restriction map
 $r^*_{l,i}:
\bH^*_{red}(T(-y_l,x_{i+1}),w)\to \bH^*_{red}(T(-y_l,x_i),w)$.
If $l\in T(-y_l,x_{i+1})\setminus T(-y_l,x_i)$ then
$l=z+E_{j(i)}$, $z\leq x_i$ and
 the coefficients of $E_{j(i)}$ in $z$ and $x_i$ are the
same. Hence, $(x_i,E_{j(i)})\geq (z,E_{j(i)})$.  This implies that
$$\chi_k(z+E_{j(i)})-\chi_k(z)\geq \chi_k(x_{i+1})-\chi_k(x_i)\geq g_{j(i)},
$$
which also shows (via \ref{RELL}(c)) that
$w_1[z,z+E_{j(i)}]=w_0(z+E_{j(i)})\geq w_0(z)$. Hence, the retract
$T(-y_l,x_{i+1})\to T(-y_l,x_i)$, which sends cycles of type
$z+E_{j(i)}$ (as above) to $z$ (and preserves all cycles of
different type) induces an isomorphism $r^*_{l,i}$. Similar
argument works if we move from $y_l$ to $y_{l+1}$. Now, property
(b) guarantees that $r^*_{-Y,X}$ is an isomorphism. On the other
hand, $\bH^*_{red}(T(-Y,X),w)$ is finitely generated over $\Z$.
\end{proof}

\para\label{CORCOR} {\bf Corollary.} {\em For any pair $(\Gamma,k)$, the space
$S_n$ is contractible for $n$ sufficiently large.}

\begin{proof}
Fix $X,Y$ as in Fact of the proof of (\ref{STR2}). Let $n$ be so large that
$T(-Y,X)\subset S_n$. Then, the same argument as in the proof of
(\ref{STR2}) shows that $S_n\cap T(-y_l,x_i)\hookrightarrow
S_n\cap T(-y_l,x_{i+1})$ admits a deformation retract. Hence, by induction,
$T(-Y,X)\subset S_n$ have the same homotopy type.
\end{proof}

If $g=0$,  one may also prove the contractibility of $S_n$ for $n\gg 0$
by verifying that $S_n$ is a deformation retract in the
real ellipsoid $ \{x\in\R^s\,:\, \chi_k(x)\leq n\}$, which is
obviously contractible.

\para\label{F} {\bf Definitions.}
We will consider the following (euler--characteristic
type) numerical invariants:
$$eu(\bH^0(\Gamma,k)):=-m_k+\rank_\Z(\bH^0_{red}(\Gamma,k)),$$
$$eu(\bH^*(\Gamma,k)):=-m_k+
\textstyle{\sum_q}(-1)^q\rank_\Z(\bH^q_{red}(\Gamma,k)).$$

\para\label{INVO} {\bf Remark.} There is a symmetry present
in the picture. Indeed, the involution $x\mapsto -x$ ($x\in L'$)
induces identities $\chi_{-k}(-l)=\chi_k(l)$, hence  isomorphisms
$$\bH^*(\Gamma,k)=\bH^*(\Gamma,-k)\ \ \mbox{and} \ \
\bH^*_{red}(\Gamma,k)=\bH^*_{red}(\Gamma,-k).$$ Notice that the
involution $[k]\mapsto [-k]$ corresponds to the natural involution
of $Spin^c_t(M)\subset Spin^c(M)$.

Regarding the canonical structure, $[K]=[-K]$ if and only if $K\in
L$. In singularity theory, such graphs are called `numerical
Gorenstein' (when the tangent bundle on $X\setminus 0$ is topologically
trivial). On the other hand, this happens if and only if the
canonical $spin^c$--structure is $spin$.

\subsection{Dependence of $\bH^*(\Gamma,k)$ on $k\in Char$.}

\para\label{3.6} Fix $\Gamma$ as above.
Above we defined for any $k\in Char$ a
graded $\Z[U]$--module  $\bH^*(\Gamma,k)$. Some of these graded
roots are not very different. Indeed,  assume that $[k]=[k']$ (cf.
\ref{2.3}), hence $k'=k+2l$ for some $l\in L$. Then
$\chi_{k'}(x-l)=\chi_k(x)-\chi_k(l)$ for any $x\in L$. Therefore,
the transformation $x\mapsto x':=x-l$ realizes the following
identification:

\para{\bf Lemma.}\label{3.7} {\em If $k'=k+2l$ for some $l\in
L$, then:} $\bH^*(\Gamma,k')=\bH^*(\Gamma,k)[-2\chi_k(l)].$

\vspace{2mm}

\noindent
 In fact,  there is an easy way to choose one module from the
multitude $\{\bH^*(\Gamma,k)\}_{k\in[k]}$. Indeed, set
$m_k=\min_{l\in L}\chi_k(l)$ as above. Since
$(k+2l)^2=k^2-8\chi_k(l)$, we get
\begin{equation*}
8m_k:=k^2-\max_{k'\in[k]}(k')^2 \leq 0.\end{equation*} Set
$M_{[k]}:=\{k\in [k]: m_k=0\}$. Hence, if
$k_0$ and $k_0+2l\in M_{[k]}$, then $-\chi_{k_0}(l)=0$. In
particular, for any fixed orbit $[k]$, any choice of $k_0\in
M_{[k]}$ provides  the same module  $\bH^*(\Gamma,k_0)$. In the
sequel we  will denote this module  by $\bH^*(\Gamma,[k])$. Notice
that with this notation, for any $k\in[k]$
\begin{equation*}
\bH^*(\Gamma,k)=\bH^*(\Gamma,[k])[2m_k].
 \end{equation*} Recall that the set of
orbits $[k]$ is the index set of the torsion $spin^c$--structures
of $M$, cf. (\ref{2.3}).

\para\label{KR1} {\bf Distinguished representative.}
There is another more sophisticated way to  choose a
representative from a class $[k]$. Let $[k]=K+2(l'+L)$. Then in
the class $l'+L$ (corresponding to an element of $H$) one can
chose $\bar{l}'_{ne}\in L'$, cf. (\ref{2.4}). The distinguished
representative of $[k]$ is, by definition,
$k_r:=K+2\bar{l}'_{ne}$. For example, the distinguished
characteristic element in $[K]$ is $K$ itself. In \cite{NOSZ}, the
elements $k_r$ had a key role. The following result basically was
proved there:

\para\label{KR2} {\bf Proposition.} {\em Fix a representative
  $k_r=K+2\bar{l}'_{ne}$ as above. Then in Fact (cf. proof of (\ref{STR2}))
one may take $Y=0$. This means that there exists an increasing sequence
$\{y_i\}_{i\geq 0}$ with $y_0=0$, $y_{i+1}=y_i+E_{j(i)}$ for some
$j(i)\in\calj$ for all $i\geq 0$, all the coefficients of $y_i$ tend to
infinity, and finally, for any $i\geq 0$ one has  }
$$\chi_{k_r}(-y_{i+1})-\chi_{k_r}(-x_i)\geq g_{j(i)}.$$

\begin{proof}
Notice that
$\chi_{k_r}(-y_{i+1})-\chi_{k_r}(-y_i)=-e_{j(i)}-1+g_{j(i)}
+(\bar{l}'_{ne}-y_i,E_{j(i)}).$

Therefore, if in a graph with $g=0$ we can find a sequence with
the wanted properties, then the same  sequence will work if we
decorate the vertices of the graph with some $g_j$. Hence, we may
assume that $g=0$. In this case the statement follows from
\cite{NOSZ} (6.1)(b), and its proof. In short, the argument is the
following. Take $Y>0$ (arbitrary large) provided by Fact. Then one
can connect $-Y$ to 0 with an increasing sequence along which
$\chi_{k_r}$ is  decreasing. Indeed, for any $y<0$ there exists
$j$ so that $E_j$ is in the support of $y$, and
$\chi_{k_r}(y+E_j)\leq \chi_{k_r}(y)$. (If not, then
$(E_j,y+\bar{l}'_{ne})\leq 0$ for all $E_j$ supported by $y$. But
the same inequality automatically works for all other components.
Hence $y+\bar{l}'_{ne}\in -L_{\Q,ne}$ with $y<0$, a
contradiction.)
\end{proof}

\para\label{KR3} {\bf Corollary.} {\em

(a) $\bH^*(\Gamma,k_r)=\bH^*((\R_{\geq 0})^s,k_r)$, i.e., in the
construction of $\bH^*(\Gamma,k_r)$ one may only work with
effective cycles from $L_e$ instead of $L$ (in other words, only
with cubes sitting in $\R_{\geq 0}^s$).

(b) With respect to the canonical characteristic element $K$,
$S_n(K)$ is connected for all $n\geq 1$.}

\begin{proof}
(a) follows from a combination of (\ref{KR2}) with the proof of
(\ref{STR2}). (b) was proved in  (6.1)(d) \cite{NOSZ} under the
assumption $c=g=0$. The very same proof (based on part (a)) can be
adopted.
\end{proof}

\subsection{(In)dependence of $\Gamma$.}

\para Clearly, many different negative definite 
plumbing graphs can provide the same
3-manifold $M$. But all these plumbing graphs can be connected by
each other by a finite sequence of blowups/downs of
$(-1)$-vertices with genus zero and  whose number of incident
edges is $\leq 2$.

\para{\bf Proposition.}\label{3.10} {\em The set
$\bH^*(\Gamma,[k])$,  where $[k]$ runs over $Spin^c_t(M)$,
depends only on $M$ and is independent of the choice of the
(negative definite) plumbing graph $\Gamma$ which provides $M$.}

\begin{proof}
First we assume that $\Gamma'$ is obtained from
$\Gamma$ by `blowing up a smooth point of one of the exceptional
curves'.  More
precisely, $\Gamma'$ denotes a graph with one more vertex and one
more edge than $\Gamma$: we glue to a vertex $j_0$ by the new
edge the new vertex with decoration $-1$ and genus 0, while the
decoration of $E_{j_0}$ is modified from $e_{j_0}$ into
$e_{j_0}-1$,  and we keep all the other decorations. We will use
the notations $L(\Gamma),\ L(\Gamma'),\ L'(\Gamma), \
L'(\Gamma')$. Similarly, write  $\, I, \ I'$ for the corresponding
intersection forms.  Set $E_{new}$ for the new base element in
$L(\Gamma')$. The following facts can be verified:

$\bullet$ Consider the maps $\pi_*:L(\Gamma')\to
L(\Gamma)$ defined by $\pi_*(\sum x_jE_j+x_{new}E_{new})=\sum
x_jE_j$, and $\pi^*:L(\Gamma)\to L(\Gamma') $  defined by
$\pi^*(\sum x_jE_j)=\sum x_jE_j +x_{j_0}E_{new}$. Then
$I'(\pi^*x,x')=I(x,\pi_*x')$.  This shows that
$I'(\pi^*x,\pi^*y)=I(x,y)$ and $I'(\pi^*x,E_{new})=0$ for any
$x,y\in L(\Gamma)$.

$\bullet$ Set the  (nonlinear)
map: $c:L'(\Gamma)\to L'(\Gamma')$,
$c(l'):=\pi^*_\Q(l')+E_{new}$. Then $c(Char(\Gamma))\subset
Char(\Gamma')$ and $c$  induces an isomorphism between the orbit
spaces $Char(\Gamma)/2L(\Gamma)$ and $ Char(\Gamma')/2L(\Gamma')$.

$\bullet$ Consider $k\in Char(\Gamma)$ and write
$k':=c(k)\in Char(\Gamma')$. Then for any $x\in L(\Gamma)$ one
has: $\chi_{k}(x)=\chi_{k'} (\pi^*x)$. 
Moreover, for any $z\in L(\Gamma')$, write $z$ in the form
$\pi^*\pi_*z+aE_{new}$ for some $a\in \Z$. Then
$\chi_{k'}(z)=\chi_{k'}(\pi^*\pi_*z)+\chi_{k'}(aE_{new})=
\chi_k(\pi_*(z))+a(a+1)/2$. Hence, the projection in the direction
$E_{new}$ provides a homotopy equivalence $S_n(\Gamma',k')\to
S_n(\Gamma,k)$.

In fact, this can be done in two steps. Let $\Pi^*$ be the union
of cubes of $\Gamma'$ with all vertices in $\pi^*(L(\Gamma))\cup
(\pi^*(L(\Gamma))-E_{new})$. Then $S_n(\Gamma',k')$ has a
deformation retract (via projection in  $E_{new}$ direction) into
$\Pi^*\cap S_n(\Gamma',k')$. On the other hand, the projection of
the later one onto $S_n(\Gamma,k)$ is a homotopy equivalence (by
checking the liftings of the cubes).

\vspace{2mm}

There is a similar verification in the case when one blows up ``an
intersection point'' corresponding to two indices $i_0$ and $j_0$
with $(E_{i_0},E_{j_0})=1$. (The only difference is that
$\pi^*(\sum x_jE_j)=\sum x_jE_j+(x_{j_0}+x_{i_0})E_{new}$.) The
details are left to the reader.
\end{proof}


\para\label{LH} {\bf Remarks.}

(a) {\bf  Lattice homology.} Obviously, there exists a parallel
homological theory as well (already used in  \cite{OSzP} for
$q=0$). Indeed, take $\calF_q:=\calC_q\otimes_\Z\et^+_0$, and
define $\partial_w:\calF_q\to \calF_{q-1}$ by
$$\partial_w(\square_q\otimes t)= \sum_k\varepsilon_k \cdot
\square_{q-1}^k\otimes U^{w(\square_q)-w(\square_{q-1}^k)}t.$$
Then $\bH_*(\calF_*,\partial_w)$ is the corresponding {\em lattice
homology} of the pair $(\Z^s,w)$. Similarly as above, it equals
$\oplus_nH_*(S_n,\Z)$. If $w$ is given as in
(\ref{dEF1}), then we get the lattice homology $\bH_*(\Gamma,k)$
of $(\Gamma,k)$.

(b) {\bf Graded root.} For each $\Gamma$ whose plumbed manifold is
rational homology sphere, and $k\in Char(\Gamma)$, the author in
\cite{NOSZ} constructed a {\em graded root}, from which one
recovers by a natural procedure $\bH^0(\Gamma,k)$. Using the
weight functions $w_0$ and $w_1$ of (\ref{DEF}), one can define in
a similar way a graded root for any $\Gamma$ (whose vertices of
degree $n$ correspond to the connected components of $S_n$) with
similar properties to those from \cite{NOSZ}.

(c) It might happen, that some non--empty real ellipsoids $\{x\in
\R^s\,:\, (\chi_k\otimes \R)(x)\leq n\}$ contain no lattice points
at all. In fact, $\min (\chi_k\otimes \R)-m_k$ can be arbitrarily
large. Take for example the rational $-2$ curve and blow up in $n$
different points. Then $\min\chi_K\otimes \R=-n/8$, but $m_K=0$.

\subsection{Path cohomology.}

\para\label{PC1} {\bf Construction.}
Fix $\Z^s$ and compatible weight functions
$w_0$ and $w_1$ as in (\ref{complex}--\ref{weight}).

We consider a sequence $\gamma:=\{x_i\}_{i=0}^t$  so that $x_0=0$,
$x_i\not=x_j$ for $i\not=j$, and $x_{i+1}=x_i\pm E_{j(i)}$ for
$0\leq i<t$. We write $T$ for the
union of 0--cubes marked by the points $\{x_i\}_i$ and
segments of type  $[x_i,x_{i+1}]$.
Then, by (\ref{SSP}), we get a graded $\Z[U]$--module $\bH^*(T,w)$,
which is called the {\em
path--cohomology} associated with the `path' $\gamma$ and weights
$\{w_q\}_q$. It is denoted by $\bH^*(\gamma,w)$.
It has an augmentation with $\et^+_{2m_\gamma}$,
where $m_\gamma:=\min_iw_0(x_i)$, and one gets the {\em reduced path
cohomology} $\bH^0_{red}(\gamma,w)$ with
$$\bH^0(\gamma,w)=\et_{2m_\gamma}^+\oplus
\bH^0_{red}(\gamma,w).$$ Similarly as in (\ref{F}), we
consider its `Euler--invariant'
$$eu(\bH^0(\gamma,w)):=-m_\gamma+\rank_\Z\bH^0_{red}(\gamma,w).$$

\para\label{PC2} {\bf Lemma.} {\em
 $\bH^q(\gamma,w)=0$ for $q\geq 1$, and }
$$eu(\bH^0(\gamma,w))=-w_0(0)+\sum_{i=0}^{t-1}\,
w_1([x_{i},x_{i+1}])-w_0(x_{i+1}).$$
\begin{proof}
Use induction comparing the paths $\{x_i\}_{i=0}^{n-1}$ and
$\{x_i\}_{i=0}^n$ ($0<n\leq t$).
\end{proof}

\para\label{PC3} {\bf Restriction map. Examples.}
In general, the restriction map $r^0:\bH^0(\R^s,w)\to
\bH^0(\gamma,w)$ is not onto. Indeed, let us fix a graph $\Gamma$
and weights as in (\ref{dEF1}), and we will study different paths
connecting $x_0=0$ with $x_t=Z_{min}$, the Artin's cycle (the
unique minimal element of $-L_{ne}\setminus 0$). In order to
simplify the picture, we assume that $c+g=0$, $k=K$, and
$\chi(Z_{min})\leq 0$ (i.e. $\Gamma$ is not rational, cf.
\ref{Ra1}). We write $\chi:=\chi_K$.  There is an `optimal' way to
find $Z_{min}$, given by Laufer's algorithm \cite{Laufer72}: start
with $x_0=0$, take for $x_1$, say, $E_1$ arbitrarily; if $x_i$
(already constructed) is in $-L_{ne}$, then stop, set $t=i$, and
$x_t=Z_{min}$; if $(x_i, E_{j(i)})>0$ for some $j(i)$ then take
$x_{i+1}=x_i+E_{j(i)}$ and continue the algorithm with $x_{i+1}$.

If one considers any path $\gamma_L$ connecting $0$ and $Z_{min}$
provided by Laufer's algorithm, then $\chi(x_1)=\chi(E_1)=1$, and
after that $\chi$ will decrease to $\chi(Z_{min})$, hence
$\bH^0(\gamma_L,K)=\et^+_{2\chi(Z_{min})}\oplus \et_0(1)$.

Assume that the multiplicity of $E_1$ in $Z_{min}$ is $\geq 2$.
Then one may take the `non--optimal' increasing path $\gamma$
connecting 0 by $Z_{min}$, by taking $x_0=0$, $x_1=E_1$,
$x_2=2E_1$, and after that we proceed according to Laufer's
algorithm. Then the maximum $\chi$--value reached is
$\chi(2E_1)=2-e_1\geq 3$ and
$\bH^0(\gamma,K)=\et^+_{2\chi(Z_{min})}\oplus \et_0(2-e_1)$.

One may verify that in the first case of $\gamma_L$ the
restriction $r^0$ is onto, while in the second case it is not.
E.g., if $\Gamma$ is minimally elliptic (see \ref{ELL}), then
$\bH^0(\Gamma,K)=\bH^0(\gamma_L,K)=\et^+_0\oplus \et_0(1)$ and
$r^0$ is an isomorphism, while in the second case $r^0$  is not
onto (by a rank argument). Moreover, in this second case,
$eu(\bH^0(\gamma,K))>eu(\bH^0(\Gamma,K))$.


\para\label{PC4} {\bf Lemma.} {\em  Fix two end--points, say 0 and
$l\in L$. We consider all the paths $\calP(l)$ as in (\ref{PC1})
with $x_0=0$ and $x_t=l$. 

(a) There exists $\gamma\in \calP(l)$ such that
$r^0:\bH^0(\Gamma,w)\to \bH^0(\gamma,w)$ is onto.

(b) If for some $\gamma\in \calP(l)$ the restriction $r^0$ is onto
then $eu(\bH^0(\gamma,w))\leq eu (\bH^0(\Gamma,w))$.\\
In particular,}
\begin{equation*}
\min_{\gamma\in \calP(l)}eu(\bH^0(\gamma,w))\leq
eu(\bH^0(\Gamma,w)).
\end{equation*}

\begin{proof} (a) Take any $\gamma$ from $\calP(l)$.
If $n\gg 0 $ then $S_n\cap\gamma$ contains all the vertices and
segments of $\Gamma$, hence it is contractible, $H^0(S_n\cap
\gamma,\Z)=\Z$ and the restriction $r^0_n: H^0(S_n,\Z)\to
H^0(S_n\cap \gamma,\Z)$ is onto. If $r^0$ is not onto, then let
$n$ be the largest integer for which $r^0_n$ is not onto. This
means that there exists $x_i$ and $x_j$ ($i<j$) so that the path
$[x_i,x_j]$ of $\gamma$ is not in $S_n$, but there is a path
$\gamma_{ij}$ connecting $x_i$ with $x_j$ in $S_n$. Then replace
$[x_i,x_j]$ by $\gamma_{ij}$. 
Notice that the higher degree homologies (for $n'>n$)  remain
unmodified. 
Repeating this procedure after a finite step we get the wished path.

The proof of (b) is left to the reader.
\end{proof}

\section{Examples.}

\noindent  The $\Z[U]$--module $\bH^0(\Gamma,k)$ is not new, it is the
combinatorial module considered in \cite{OSzP,NOSZ} (where
it was denoted by $\bH^+$). In fact, regarding $\bH^0(\Gamma,k)$,
\cite{NOSZ} is one of our main sources of examples.

\subsection{The case of rational graphs.}

\para\label{Ra1} {\bf Definition.} By its very definition, a singularity is
rational of its geometric genus $p_g$ is vanishing. By
\cite{Artin62,Artin66}, this can be characterized combinatorially:
a singularity is rational if and only if its graph satisfies
$\min_{l\in L_e\setminus 0}\chi_K(l)>0$ (or, equivalently,
$\chi_K(Z_{min})=1$, or $\chi_K(Z_{min})>0$). Therefore,  a
(connected, negative definite) graph with this property is
 called rational.  For them  $c=g=0$ automatically.

In \cite{NOSZ} the author have given a different characterization:

\para\label{Ra2} {\bf Proposition.} {\em  Assume that $\Gamma$ is
a negative definite connected plumbing graph with $c=g=0$. Then
$\Gamma$ is rational if and only if \, $\bH^0_{red}(\Gamma,K)=0$.
Moreover, in this case, $\bH^0_{red}(\Gamma,k)=0$, and also
$m_k=0$,  for any $k\in Char$. }

\vspace{2mm}

Even if one drops the assumption $c=g=0$, one can prove:

\para\label{Ra2b} {\bf Proposition.} {\em If \,
$\bH^0(\Gamma,K)=\et^+_0$ then $\Gamma$ is rational.}

\begin{proof}
Using (\ref{KR3}) we get $\chi_K|L_e\geq 0$. Hence
$\chi_K(E_j)=1-g_j\geq 0$. Assume that $\chi_K(Z_{min})=0$. Then
$Z_{min}$ cannot be connected by 0 in $S_0\cap L_e$ since
$w_1([0,E_j])=1$ by (\ref{RELL})(b). This contradicts the
assumption, hence $\chi_K(Z_{min})>0$, i.e. $\Gamma$ is rational.
\end{proof}

We add to this the following vanishing result:

\para\label{Ra3} {\bf Proposition.} {\em If \, $\Gamma$ is rational then
 $\bH^*_{red}(\Gamma,k)=0$ for any $k\in Char$.}

\begin{proof}
By (\ref{3.7}) we may replace $k$ with any characteristic element
in its class. Take the distinguished representative $k_r$, cf.
(\ref{KR1}).  The result follows from the proof of (\ref{STR2})
once we show that one may take $X=Y=0$ in Fact. By (\ref{KR2}) we
may take $Y=0$. Hence, we have to show that for $\Gamma$ rational
there exists an increasing sequence $\{x_i\}_{i\geq 0}$ with
$x_0=0$, $x_{i+1}=x_i+E_{j(i)}$, all the coefficients of $x_i$
tend to infinity, and $\chi_{k_r}(x_{i+1})\geq \chi_{k_r}(x_i)$.
For this take a sequence $\{z_i\}_{i=0}^t$ which connects 0 and
$Z_{min}$ provided by Laufer's algorithm, cf. (\ref{PC3}). Then
$z_0=0$, $z_1=E_1$, and $(E_{j(i)},z_i)=1$ for $1\leq i<t$
\cite{Laufer72}. Hence $\chi_{k_r}(z_1)=1-(\bar{l}'_{ne},E_1)\geq
1$ and
$\chi_{k_r}(z_{i+1})-\chi_{k_r}(z_i)=-(\bar{l}'_{ne},E_{j(i)})\geq
0$. Hence, the sequence $\{mZ_{min}+z_i\}$ with $m\geq 0$ and
$0\leq i < t$ works.
\end{proof}

Therefore, the above proof combined with the proof of
(\ref{CORCOR}) gives

\para\label{Ra4} {\bf Corollary.} {\em If \, $\Gamma$ is rational,
then $S_n$ is contractible for any $k\in Char$ whenever is
non--empty. }

\subsection{The case of elliptic graphs.}\label{ELL}

\para\label{El1} {\bf Definition.} \cite{Laufer77,Wagreich}
 A connected negative definite
graph is elliptic if $\min_{l\in L_e\setminus 0}\chi_K(l)=0$.

\vspace{2mm}

In this case $\Gamma$ might have a cycle or a vertex with genus
one, but in any case $c+g\leq 1$. The next characterization result
was proved in \cite{NOSZ} for $c=g=0$, here we verify the general
situation.

\para\label{El2} {\bf Proposition.} {\em $\Gamma$ is elliptic if and only if
\, $\bH^0(\Gamma,K)=\et_0^+\oplus \et_0(1)^\ell$ for some $\ell\geq
1$. }

\begin{proof}
Notice that, by (\ref{KR3}), $\bH^0(\Gamma,K)=\bH^0(\R_{\geq
0}^s,K)$. Assume that $\Gamma$ is elliptic. Then $\chi_K|L_e\geq
0$. Moreover, by (\ref{KR3}), $S_n(K)$ is connected for $n\geq 1$, hence
$\bH^0(\Gamma,K)=\et^+_0\oplus \et_0(1)^\ell$ for some $\ell\geq
0$. Since $\Gamma$ is not rational, $\ell\not=0$. Conversely, if
$\bH^0$ has that form, then $\chi_K|L_e\geq 0$ and there exists a
cycle $x\in L_e\setminus 0$ with $\chi_K(x)=0$, hence $\Gamma$ is
elliptic.
\end{proof}

 In fact, in the `classical' theory of elliptic singularities, there is a
combinatorial integer which guides the main topological and
analytical  properties, namely, {\em  the length of the elliptic
sequence} $\ell^{es}$, introduced by Laufer and S. S.-T. Yau (see
e.g. \cite{Yau4,Yau1}). E.g., Yau proved that $\ell^{es}+1$ is a
topological upper bound for the geometric genus, and \cite{Ninv}
shows that it  is realized by any Gorenstein singularity when
$c=g=0$.   The `simplest' elliptic singularities, the minimally
elliptic ones, are characterized by $\ell^{es}=0$, or by the
identity $Z_{min}=-K$  \cite{Laufer77}.

The point is that the above integer $\ell$ provided by
(\ref{El2}), in fact, equals  $\ell^{es}+1$. In particular, for
minimally elliptic singularities one has $\ell=1$.

We exemplify the above proposition (\ref{El2}) for three minimally
elliptic singularity. The simplest case, when $c=g=0$,  is the
hypersurface singularity $\{x^2+y^3+z^7=0\}$, whose minimal good
resolution graph has four vertices and three edges with
$E_1^2=-1,\ E_2^2=-2,\ E_3^2=-3,\ E_4^2=-7$, and  $E_1$ is
connected with the others. Then $\ker(U)$ has rank 2,  the
generating lattice points are the zero cycle and
$Z_{min}=6E_1+3E_2+2E_3+E_4$. (It general, the generators of
$\ker(U)$ correspond to  lattice points important in singularity
theory as well, cf. \cite{NOSZ,[65]}.)

\para\label{El4} {\bf Example.} Assume that $\Gamma$ consists of
three vertices, each pair of vertices is connected by an edge, the
self--intersections are $-2,\ -2, \ -3$, and $g=0$. 
In this case $-K=\sum E_j=Z_{min}$, hence $\Gamma$ is
minimally elliptic \cite{Laufer77}. With the
notation $l=\sum_jx_jE_j$, $2\chi_K(l)=
(x_1-x_2)^2+(x_2-x_3)^2+(x_3-x_1)^2+x_3^2-x_3$,
which  is always non--negative on $L$. It is zero at $(0,0,0)$ and $(1,1,1)$.
They cannot be connected by a segment, hence $S_0$ consists of
two points. On the other hand, there are a lot of lattice points
$l$ in $S_1$: with the third coordinate $x_3=-1$ one has the
pair $(x_1,x_2)=(-1,-1)$, with $x_3=0$  the pairs
$P=\{(-1,-1),\ (-1,0),\ (0,-1),\ (0,0),\
(1,0),\ (0,1),\ (1,1)\}$; with  $x_3=1$  the
pairs $P+(1,1)$, and for $x_3=2$ the pair $(2,2)$. (They are
situated symmetrically with respect to $-K/2$.) Hence, $S_1$
consists of 16 lattice points. One can verify that they can be
connected by segments, $S_1$ contains eight 2--cubes and one
3--cube, and  $S_1$ is contractible. And this is the case
for all $S_n$ with $n\geq 1$.

\para\label{El5} {\bf Example.} Assume that $\Gamma$ has only one
vertex with self--intersection $-1$ and $g=1$. In this case again
$K=-E$ and the graph is minimally elliptic (it is
the `simple--elliptic' singularity $\tilde{E}_8$:
$\{x^2+y^3+z^6=0\}$). Notice that there are only two lattice
points $l$ with $\chi_K(l)=0$, namely the zero cycle and $E$, and
they are connected by the segment $[0,E]$ from $\calQ_1$. Hence,
$w_1[0,E]=0$ would imply  $\bH^0_{red}(\Gamma,K)=0$. Therefore, in
order to have the `right' result, we are forced to put 1 for the
weight of this segment, a fact compatible with
(\ref{dEF1})-(\ref{RELL}).

Then, with this weight functions, one has: $\bH^q=0$ for $q>0$ and
$\bH^0(\Gamma,K)=\et_0^+\oplus \et_0(1)$.

\vspace{2mm}

Notice that in the case of a  minimally elliptic singularity, $K$
is integral \cite{Laufer77}, hence $[K]=[-K]$  (cf. \ref{INVO}).
We may add the following vanishing result for the other
$spin^c$--structures:

\para\label{El6} {\bf Proposition.} {\em If \, $\Gamma$ is
minimally elliptic, and the associated minimal (resolution) graph
is good,  then $\bH^*_{red}(\Gamma,[k])=0$ for any $[k]\not=[K]$.}

\begin{proof} By (\ref{3.10}) we may assume that $\Gamma$ is
minimal. Then the proof of (\ref{Ra3}) can be adopted. Indeed,
consider the representative $k_r$. Since $[k_r]\not=[K]$,
$k_r\not=0$, hence there exists at least one $j$ with
$(\bar{l}'_{ne},E_j)<0$. By \cite{Laufer77} (p. 1261-1262), there
exists a computation sequence $\{z_i\}_{i=0}^t$ for $Z_{min}$ so
that the last $E_{j(t-1)}$ is $E_j$, $(z_i,E_{j(i)})=1$ for $1\leq
i\leq t-2$, and $(z_{t-1},E_{j(t-1)})=2$ (this fact uses the
minimality of $\Gamma$). Then the proof of (\ref{Ra3}) works in
this case too.
\end{proof}

\subsection{The case of almost rational graphs.}

\para\label{AR1} {\bf Definition.} \cite{NOSZ} Assume that the graph
$\Gamma$ is connected and  negative definite with $c+g=0$. We say
that $\Gamma$ is {\em almost-rational} if there exists a vertex
$j_0\in\calj$ of $\Gamma$ such that replacing its Euler number
$e_{j_0}$ by some $e_{j_0}'\leq e_{j_0}$
 we get a rational graph.
 (In general, the choice of $j_0$ is not unique.)

\para\label{AR2}{\bf Examples.} Almost rational graphs include:
rational graphs, elliptic graphs (with $c+g=0$),  star-shaped
graphs (with central vertex of genus zero). 
But there are more `exotic' ones as well; e.g. the plumbing graph
of the rational surgery 3--manifolds $S^3_r(K)$, where
$r\in\Q_{<0}$ and $K$ is an algebraic knot in $S^3$ (see e.g.
\cite{[65],trieste}). On the other hand, not every graph is almost
rational. For example, if $\Gamma$ has two (or more) vertices $j$
with $-e_j+2$ less than or equal to the valency of the vertex $j$,
then $\Gamma$ is not almost rational (e.g. the graph from
(\ref{NV1})).

For almost rational graphs,  $\bH^0$ might be rather complicated module
(see e.g. \cite{NOSZ} for the
explicit description in the case of star--shaped graphs).
 On the other hand, we have:

\para\label{AR3} {\bf Theorem.} {\em For any almost rational
graph, $\bH^q(\Gamma,k)=0$ for any $q>0$ and $k\in Char$.}

\begin{proof}
The complete proof is rather technical and long, and we will omit
it. It is based on the results of \cite{NOSZ}, \S 9. In fact, one
can construct an infinite increasing path $\gamma$:
$\{x_i\}_{i\geq 0}$ with $x_0=0$, so that the restriction
$r^*:\bH^*(\Gamma,k)\to \bH^*(\gamma,k)$ is an isomorphism (the
fact that $r^0$ is an isomorphism is the main result of \cite{NOSZ}
\S 9). The isomorphism
 is induced by a deformation retract whose existence is proved
by a combination of results from \cite{NOSZ} with the proof of
(\ref{Ra3}).
\end{proof}

\subsection{Examples with non--vanishing $\bH^1$.}

\para\label{NV1} {\bf Example.}
Consider the following graph:

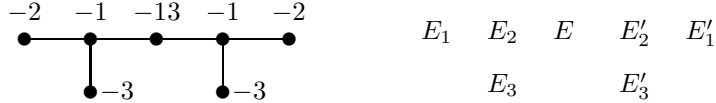
\begin{picture}(300,45)(70,0)
\put(125,25){\circle*{5}} \put(150,25){\circle*{5}}
\put(175,25){\circle*{5}} \put(200,25){\circle*{5}}
\put(225,25){\circle*{5}} \put(150,5){\circle*{5}}
\put(200,5){\circle*{5}} \put(125,25){\line(1,0){100}}
\put(150,25){\line(0,-1){20}} \put(200,25){\line(0,-1){20}}
\put(125,35){\makebox(0,0){$-2$}}
\put(150,35){\makebox(0,0){$-1$}}
\put(175,35){\makebox(0,0){$-13$}}
\put(200,35){\makebox(0,0){$-1$}}
\put(225,35){\makebox(0,0){$-2$}} \put(160,5){\makebox(0,0){$-3$}}
\put(210,5){\makebox(0,0){$-3$}}

\put(275,25){$E_1$} \put(300,25){$E_2$} \put(325,25){$E$}
\put(350,25){$E_2'$} \put(375,25){$E_1'$} \put(300,5){$E_3$}
\put(350,5){$E_3'$}
\end{picture}

On the right hand side we give names to the base elements. Set
$\chi:=\chi_K$.  We prefer to write any $l\in L$ in the
form $l=l_x+ zE+ l_y$, where $l_x=\sum x_iE_i$, $l_y=\sum
y_iE_i'$, $x_i,\ y_i,\ z\in \Z$ ($i=1,2,3$); or in the form
$(x_1,x_2,x_3;z;y_1,y_2,y_3)$. Then $-K=(7,14,5;3;7,14,5)$ and
$Z_{min}=(3,6,2;1;3,6,2)$, with
$\chi(Z_{min})=\chi(2Z_{min})=-1$. In fact, $m_K=-1$ too.
Then, it turns out that
$$\bH^0(\Gamma,K)=\et_{-2}^+\oplus \et_{-2}(1)\oplus
\et_0(1)\oplus \et_0(1),$$ where the generators of $\ker(U)$ with
homogeneous degree $-2$ are (the dual classes of) $Z_{min}$ and
$2Z_{min}$, while with degree 0 are the (dual classes of) zero
cycle and $-K$. Moreover, there exists a non--trivial class in
$\bH^1$ of homogeneous degree 0. In fact,
$$\bH^1(\Gamma,K)=\et_0(1),\ \ \mbox{and} \ \ \bH^q(\Gamma,K)=0\ \ \mbox{for
$q\geq 2$.}$$

In order to see (at least part of) these, we will analyse $S_{-1}$
and $S_0$. Since $\chi(l)=\chi(-K-l)$, we can use all the time the
$\chi$--symmetry of the lattice points with respect to $-K/2$.
If $z=0$ then
$\chi(l)=\chi(l_x)+\chi(l_y)$, and since $l_x$ and $l_y$ are
supported by rational subgraphs, $\chi(l_x)\geq 0$, $\chi(l_y)\geq
0$. Hence, the lattice points in  $S_{-1}$
have $z=1$ or $z=2$, and they correspond by the above symmetry.

Let us assume that $z=1$. Then
$\chi(l)=\chi(l_x)+\chi(l_y)-x_2-y_2+1$. Therefore, with the
notation $$f(x):=\chi(l_x)-x_2=(\
2x_1^2+x_2^2+3x_3^2-2x_1x_2-2x_2x_3-x_2-x_3\,)/2,$$ we have
$\chi(l)=f(x)+f(y)+1$. By real calculus the minimum of $f$ over
$\R^3$ is $>-2$, hence its minimum over $\Z^3$ is $\geq -1$.
Therefore, $\chi(l)=-1$ if and only if $f(x)=f(y)=-1$. By a
computation, the integral solutions of $f(x)=-1$ are the triplets
$$A:=\{(1,2,1),(1,3,1),(2,3,1),(2,4,1),(2,4,2),(2,5,2),(3,5,2),(3,6,2)\}.$$
Therefore, points $(x,1,y)$ with $x\in A$ and $y\in A$ (denoted
simply by $(A,1,A)$) are in $S_{-1}$. Let $B=(7,14,5)-A$. Then, by
symmetry, we get that the set of lattice points of $S_{-1}$ is
$(A,1,A)\cup(B,2,B)$. They determine two contractible
connected components of $S_{-1}$ in which $Z_{min}$ and $2Z_{min}$
are `representatives'.

Next, we plan to solve the equation $\chi(l)=0$ with $z=1$. Then,
$f(x)+f(y)=-1$. $f(x)=0$ has 24 integral solutions, namely union
of the triplets $A':=$
$$\{(0,0,0),(0,1,0),(1,1,0),(1,2,0),(0,1,1),(1,1,1),
(0,2,1),(2,2,1), (1,4,1),(3,4,1),(2,5,1),(3,5,1)\},$$ and the
triplets of type $A''=(4,8,3)-A'$. Set $\tilde{A}:=A\cup A'\cup
A''$, and $\tilde{B}=(7,14,5)-\tilde{A}$. Then the points of type
$$X:=(A,1,\tilde{A})\cup (\tilde{A},1,A)\cup (B,2,\tilde{B})\cup
(\tilde{B},2,B)$$ are in $S_0$. Since  $\tilde{A}\cap
B$ is not empty, all the points from $X$ can be connected by
segments. In fact, $S_0$ has three connected components, one of
them contains the zero cycle, the other contains $-K$, and the
third one, $CS_0$, contains all the points from $X$.

Finally, notice that the two intersection points $P=\tilde{A}\cap
B=(4,8,3)$ and $Q=A\cap \tilde{B}=(3,6,2)$ create a loop in
$CS_0$. Indeed, half of it is the connecting path of $(P;1;Q)$ and
$(Q;1;P)$ through points in $X$ with $z=1$, the other half 
connects  $(P;2;Q)$ with $(Q;2;P)$ through points in $X$ with
$z=2$. This loop can be contracted only in $S_1$ (which is
contractible).

\para\label{NV2} {\bf Example.} In the above example
the subgraph of $\{E_1,E_2,E_3\}$ is a `cusp', $\Gamma$ was
obtained by gluing two cusps to the `central' curve $E$. One may
create non--trivial higher dimensional modules by gluing $k$ cusps
to a central curve $E$ which has self--intersection $-6k-1$.

\section{Heegaard--Floer homology and singularity links.}

\subsection{Heegaard--Floer homology.}\

\vspace{1mm}

 \noindent  In this section we will assume that $M$ is
an oriented rational homology 3--sphere.

\para\label{HF1} {\bf Review.} Heegaard--Floer homology $HF^+(M)$
was introduced by Ozsv\'ath and Szab\'o in \cite{OSz} (and
intensively studied in a series  of articles). $HF^+(M)$ is a
$\Z[U]$-module with a $\Q$-grading compatible with the
$\Z[U]$-action, where $\deg(U)=-2$. Additionally, $HF^+(M)$ also
has an (absolute) $\Z_2$-grading; $HF^+_{even}(M)$, respectively
$HF^+_{odd}(M)$, denote the part of $HF^+(M)$ with the
corresponding parity. Moreover, $HF^+(M)$  has a natural direct
sum decomposition of $\Z[U]$-modules (compatible with all the
gradings) corresponding to the $spin^c$-structures of $M$:
$$HF^+(M)=\oplus_{\sigma\in Spin^c(M)}\ HF^+(M,\sigma).$$
For any $spin^c$-structure $\sigma$, one has a graded
$\Z[U]$-module isomorphism
$$HF^+(M,\sigma)=\calt^+_{d(M,\sigma)}\oplus HF^+_{red}(M,\sigma),$$
where 
$HF^+_{red}(M,\sigma)$ has a
finite $\Z$-rank and an induced (absolute) $\Z_2$-grading. One
also considers
$$\chi(HF^+(M,\sigma)):=\rank_\Z HF^+_{red,even}(M,\sigma)
-\rank_\Z HF^+_{red,odd}(M,\sigma).$$ Then one recovers the
Seiberg-Witten topological invariant of $(M,\sigma)$ (see
\cite{Rus}) via
$$\ssw(M,\sigma):=\chi(HF^+(M,\sigma))-d(M,\sigma)/2.$$
With respect to the change of orientation the above invariants
behave as follows: The $spin^c$-structures $Spin^c(M)$ and
$Spin^c(-M)$ are canonically identified (where $-M$ denotes $M$
with the opposite orientation). Moreover,
$d(M,\sigma)=-d(-M,\sigma)$ and
$\chi(HF^+(M,\sigma))=-\chi(HF^+(-M,\sigma))$. Notice also that
one can recover $HF^+(M,\sigma)$ from $HF^+(-M,\sigma)$ via (7.3)
\cite{OSz} and (1.1) \cite{OSzTr}.

\para\label{HF2} {\bf Example.} If $M$ is an integral homology sphere
then for the unique (=canonical) $spin^c$--structure
$\sigma_{can}$,   $\ssw(M,\sigma_{can})$ equals the Casson
invariant $\lambda(M)$ (normalized as in \cite{Lescop} (4.7)).

\subsection{Lattice homology and Heegaard--Floer homology.}\label{ARHF}

\para Assume that $\Gamma$ is a connected negative definite
plumbing graph whose associated plumbed 3--manifold is a {\em rational
homology sphere}. Our goal is to recover the Heegaard--Floer
homology of $M$ in a purely combinatorial way from $\Gamma$.
We write $\#\calj=s$.

\para\label{HF3}{\bf Theorem.} \cite{OSzP,NOSZ} {\em Assume that ~$\Gamma$ is an
almost rational graph. Then, for any $[k]\in Spin^c(M)$
$$HF^+_{odd}(-M,[k])=0,$$ and
$$HF^+_{even}(-M,[k])=\bH^0(\Gamma,[k])\Big[-\max_{k'\in [k]}\frac{(k')^2+s}{4}\Big].$$
In particular (cf. \ref{3.7}), for any $k\in[k]$ one has}
\begin{equation*}
d(M,[k])=\max_{k'\in[k]}\frac{(k')^2+s}{4}=\frac{k^2+s}{4}-2\min\chi_k.
\tag{$*$}\end{equation*}

\para\label{HF4}{\bf Corollary.} {\em If \, $\Gamma$ is an almost rational graph,
then for any $k\in Char$:}
$$-\ssw(M,[k])-\frac{k^2+s}{8}=-\min\chi_{k}+\rank_{\Z}\bH_{red}^0(\Gamma,k)
=eu(\bH^0(\Gamma,k)).$$
\para\label{HF5} {\bf Conjecture.} Let $M$ be a plumbed rational homology sphere
associated with a connected negative definite graph $\Gamma$. Then
for any $k\in Char$ the identity ($*$) of (\ref{HF3}) is valid, and 
\begin{equation*}
-\ssw(M,[k])-\frac{k^2+s}{8}=-\min\chi_{k}+\sum_q(-1)^q
\rank_{\Z}\bH_{red}^q(\Gamma,k)=eu(\bH^*(\Gamma,k)).\tag{$**$}
\end{equation*}
\noindent
 In fact, we predict that with $d=d(M,[k])$:
 $$HF^+_{red,even}(-M,[k])=\bigoplus_{q
\ even} \bH^q_{red}(\Gamma,[k])[-d], \ \ \mbox{and} \ \ \
HF^+_{red,odd}(M,[k])=\bigoplus_{q\ odd} \bH^q_{red}(\Gamma,[k])[-d].$$

\para\label{HF6} {\bf Example.} Take $\Gamma$ from (\ref{NV1}), and $k=K$.
Then ($*$) is true by \cite{OSzP}, Corollary 1.5. Moreover, by
(\ref{NV1}), $eu(\bH^*)=-(-1)+3-1=3$, $(K^2+s)/8=-1$. On the other hand,
the Casson invariant of $M$ is $-2$ (using, e.g.,  the formula of Ra\c{t}iu,
see (5.3)  \cite{[51]}). Hence, ($**$) is valid as well.

\para\label{HF7} {\bf Remarks.} (a) The above identities are not valid (in this
form) when $c+g>0$. 

(b)   (\ref{Ra2}--\ref{Ra3}), or (\ref{HF3}) shows that if $\Gamma$ is
rational then  $M$ is an $L$--space (in the sense of
Ozsv\'ath and Szab\'o, i.e. $HF^+_{red}(M)=0$).
From the perspective of Conjecture (\ref{HF5}), we expect that this is an
`if and only if' correspondence: $\Gamma$ is rational if and only if $M$
is an $L$--space. Notice that by (\ref{Ra2}), if $\bH°_{red}=0$ then
$\Gamma$ is rational.

(c) Although we expect an identification of the $\bH^*$ modules with
the Heegaard--Floer modules $HF^+$, the lattice cohomology (apparently)
contains more structure (at least, the author is not able to recover them
in $HF^+$). For their existence the explanation is, maybe, that the involved
3--manifolds are rather special.
We list here three such extra properties.

(i) The (absolute) grading of $\{\bH^q\}_{q\geq 0}$ (with respect to $q$)
is indexed by $\Z$ in contrast with the $\Z_2$ (even/odd) grading of $HF^+$.

(ii) Consider from (\ref{STR1}) the identity $\bH^*=\oplus_nH^*(S_n,\Z)$.
How can the ring structure of each $H^*(S_n,\Z)$ exploited?

(iii) For a fixed graph $\Gamma$,
consider any distinguished representative $k_r$. Since
$\chi_{k_r}(l)\geq \chi_K(l) $ for any $l\in L$, we get
$S_n(k_r)\subset S_n(K)$, hence a  natural ring homomorphism
$H^*(S_n(K))\to H^*(S_n(k_r))$, or $R: \bH^*(\Gamma,K)\to
\bH^*(\Gamma,k_r)$. Notice that in the case of rational or elliptic graphs it
happens that properties of the  module associated with the canonical
$spin^c$--structure `dominates' all the others, cf. (\ref{Ra2}) or
(\ref{El6}). Is it possible to say something similar in general? Is $R$ onto?

\section{Line bundles associated
with surface singularities.}\label{linebundles}

\noindent Starting from this section, we start to analyse the
analytic aspects of the singularity $(X,0)$  as well. The analytic type is
preserved in the complex manifold structure of the resolution
$\tx$. Holomorphic line bundles on $\tx$ codify a lot of
information about it.

\subsection{Cohomological computations}

\para\label{CohCom} Let \ $\pi:(\tx,E)\to (X,0)$ be a fixed good resolution of
$(X,0)$. Let $Pic(\tx)$ be the group of isomorphism classes of
holomorphic line bundles on $\tx$  and  $c_1:Pic(\tx)\to L'$,
$c_1({\call})=\sum_j \deg(\call|E_j)\, D_j$ the
set of Chern classes  of $\call$.
We prefer to use the same notation for $l=\sum n_jE_j\in L $ and
divisors  $\sum n_jE_j $ of $\tx$ supported by $E$. Hence,
we can consider the line bundle
$\calo_{\tx}(l):=\calo_{\tx}(\sum n_jE_j)$. If $l>0$, we write
$\chi(l)$ for $\chi_K(l)=\chi(\calo_l)$ (cf. (\ref{chifunction})).
We write $|l|$ for the support of $l$.

In this subsection we analyse $h^1(\call):=\dim H^1(\tx,\call)$
for any $\call\in Pic(\tx)$.
First, recall  the following general (Grauert-Riemenschneider
type) vanishing theorem (cf.\ \cite{MR}, page 119, Ex. 15):

\para\label{CC1}
 {\em If $c_1(\call)\in K+L_{\Q,ne}$, then $h^1(l,\call|_l)=0$ for
any $l\in L$, $l>0$, hence $h^1(\tx,\call)=0$. }

\vspace{1mm}

\noindent The next statement is an improvement of it, valid for
rational singularities:

\para\label{4.3} {\em Assume that $(X,0)$ is a rational singularity. If
$c_1(\call)\in L_{\Q,ne}$, then $h^1(l,\call|_l)=0$ for any $l>0$,
$l\in L$, hence $h^1(\tx,\call)=0$ too. } \
\begin{proof} From the point of view of the next discussion, it is instructive to
  see the proof.  For any $l>0$ there
exists $E_j\subset |l|$ such that $(E_j,l+K)<0$. Indeed,
$(E_j,l+K)\geq 0$ for any $j$ would imply $\chi(l)=-(l,l+K)/2\leq
0$, which would contradict the rationality of $(X,0)$ \cite{Artin62}.
Then,  using
$$0\to \call\otimes \calo_{E_j}(-l+E_j)\to \call|_l\to \call|_{l-E_j} \to 0$$
one gets $h^1(\call|_l)=h^1(\call|_{l-E_j})$, hence by induction
$h^1(\call|_l)=0$.
\end{proof}
This will be generalized in two different ways. First we show that
the computation of any $h^1(\call)$ can be reduced  to the
computation of  some $h^1(\call')$ with $c_1(\call')\in
L_{\Q,ne}$.

\para\label{4.4}{\bf Proposition.} {\em Let $\tx\to X$ be a good
resolution of a normal singularity $(X,0)$ as above.

(a) For any $l'\in L'$ there exists a unique minimal element
$l_{l'}\in L_e$ with $e(l'):=l'-l_{l'}\in L_{\Q,ne}$.

(b) $l_{l'}$ can be found by the following (generalized Laufer's)
algorithm. One constructs a sequence  $x_0, x_1, \ldots, x_t\in
L_e$ with $x_0=0$ and $x_{i+1}=x_i+E_{j(i)}$, where each  index
$j(i)$ is determined by the following principle. Assume that $x_i$
is already constructed. Then, if $l'-x_i\in L_{\Q,ne}$, then one
stops, and $t=i$. Otherwise, there exists at least one $j$ with
$(l'-x_i,E_j)<0$. Take for $j(i)$ one of these $j$'s. Then this
algorithm stops after a finitely many steps, and $x_t=l_{l'}$.

(c) For any $\call\in Pic(\tx)$ with $c_1(\call)=l'$ one has:
$$h^1(\call)=h^1(\call\otimes \calo_{\tx}(-l_{l'}))-(l',l_{l'})-\chi(l_{l'}).$$
In particular (since $c_1(\call\otimes \calo_{\tx}(-l_{l'}))\in
L_{\Q,ne}$), the computation of any $h^1(\call)$ can be reduced
(modulo the combinatorics of $L$)
to the computation of  some $h^1(\call')$ with $c_1(\call')\in
L_{\Q,ne}$. }

%
\begin{proof} (a) Since $(\,,)$ is negative definite,
there exists $l\in L_e$  with $l'-l\in
L_{\Q,ne}$ (take e.g. a large multiple of some $Z$ with
$(Z,E_j)<0$ for any $j$). Next, we  prove that if $l'-l_i\in
L_{\Q,ne}$ for $l_i\in L_e$, $i=1,2$, and $l:=\min\{l_1,l_2\}$,
then $l'-l\in L_{\Q,ne}$ as well. For this, write $x_i:=l_i-l\in
L_e$. Then $|x_1|\cap|x_2|=\emptyset$, hence for any fixed $j$,
$E_j\not\subset |x_i|$ for at least one of the $i$'s. Therefore,
$(l'-l,E_j)=(l'-l_i,E_j)+(x_i,E_j)\geq 0$.

(b) First we prove that $x_i\leq l_{l'}$ for any $i$. For $i=0$
this is clear. Assume that it is true for some $i$ but not for
$i+1$, i.e.\ $E_{j(i)}\not\subset |l_{l'}-x_i|$. But this would
imply
$(l'-x_i,E_{j(i)})=(l'-l_{l'},E_{j(i)})+(l_{l'}-x_i,E_{j(i)})\geq0$,
a contradiction. The fact that $x_i\leq l_{l'}$ for any $i$
implies that the algorithm must stop, and $x_t\leq l_{l'}$. But
then by the minimality of $l_{l'}$ (part a) $x_t=l_{l'}$. (Cf.\
\cite{Laufer72}.)

(c) For any $0\leq i<t$, consider the exact sequence
$$0\to \call\otimes \calo_{\tx}(-x_{i+1}) \to
\call\otimes \calo_{\tx}(-x_{i}) \to \call\otimes
\calo_{E_{j(i)}}(-x_{i}) \to 0.$$ Since $\deg(\call\otimes
\calo_{E_{j(i)}}(-x_{i}))=(l'-x_i,E_{j(i)})<0$, one gets
$h^0(\call\otimes \calo_{E_{j(i)}}(-x_{i}))=0$. Therefore
$$h^1(\call\otimes \calo_{\tx}(-x_{i}))-
h^1(\call\otimes \calo_{\tx}(-x_{i+1}))= -\chi(\call\otimes
\calo_{E_{j(i)}}(-x_{i}))$$ which equals
$-(l',x_{i+1}-x_i)+\chi(x_i)-\chi(x_{i+1})$. Hence the result
follows by induction.
\end{proof}

\para\label{4.5e}{\bf Examples. Rational singularities.}  If $(X,0)$
is rational then $c_1:Pic(\tx)\to
L'$ is an  isomorphism. Moreover, using (\ref{4.3}) and (\ref{4.4})(c), one has
$h^1(\call)=-(l',l_{l'})-\chi(l_{l'})$.
 In particular, $h^1(\call)$ {\em depends only
on $\Gamma$ and it is independent of the analytic structure of
$(X,0)$}.

\subsection{Path cohomology and upper bounds for
$h^1(\call)$.}\label{PCUB}

\para\label{UJ1} For the next result, we start with the following
set of data and notations: $\call\in Pic(\tx)$, $l':=c_1(\call)$,
$k:=K-2l'$ (cf. \ref{chifunction}). We consider  a `path'
$\gamma$: $\{x_i\}_{i=0}^t$, where
$x_0=0$, $x_t\in l'-K-L_{\Q,ne}$, and $x_{i+1}=x_i\pm E_{j(i)}$
for some $j(i)\in\calj$ ($0\leq i<t$).

Using the exact sequence $0\to \call\otimes \calo(-x_t)\to
\call\to \call|_{x_t}\to 0$, and the Grauert-Riemenschneider
vanishing (\ref{CC1}), we get  $h^1(\call)=h^1(\call|_{x_t})$
(this motivates the corresponding restriction for $x_t$). In the
next proposition the `symbol' $h^1(\call|_{x_0})$ will stand  for
zero.

\para\label{UJ2} {\bf Proposition.}  {\em With the above
notations, for any $0\leq i<t$  with $x_{i+1}>x_i$ one has:
$$h^1(\call|_{x_{i+1}})-h^1(\call|_{x_i})\leq \left\{\begin{array}{ll}
-\Delta_i+M^{g_{j(i)}}(-\Delta_i) & \mbox{if \ $\Delta_i<0$},\\
M^{g_{j(i)}}(\Delta_i)& \mbox{if \ $\Delta_i\geq 0$},
\end{array}\right. $$where
$\Delta_i:=\chi_k(x_{i+1})-\chi_k(x_i)$. If $x_{i+1}<x_i$ then
$h^1(\call|_{x_{i+1}})-h^1(\call|_{x_i})\leq 0$. }

\vspace{2mm}

In particular, adding all these inequalities, we get a topological 
upper bound 
for $h^1(\call)$.

\para\label{UJ3} {\bf Example.} {\em Assume that $g=0$ and
$\gamma$ is increasing.
Since $M^0(n)=0$ for all $n\geq 0$, we get}
$$h^1(\call)\leq \sum_{i=0}^{t-1} \max\{\, 0\, ,\,
\chi_k(x_i)-\chi_k(x_{i+1})\,\}.$$

\noindent {\em Proof of (\ref{UJ2}).} 
Assume that $x_{i+1}>x_i$ (the other case is trivial). Write
$\calm_i$ for the line bundle $\call\otimes
\calo_{E_{j(i)}}(-x_i)$ on $E_{j(i)}$. From the cohomological
exact sequence
$$\cdots \to H^1(E_{j(i)},\calm_i)\to
H^1(\call|_{x_{i+1}})\to H^1(\call|_{x_i})\to 0$$ we have to
estimate $h^1(\calm_i)$. Notice that $\chi(\calm_i)=\Delta_i$
by (\ref{chifunction}). Hence, if $\Delta_i\geq 0$, then
$h^1(\calm_i)\leq M^{g_{j(i)}}(\Delta_i)$, by the very definition
of $M^g(n)$. Assume that $\Delta_i<0$. Then, by Serre duality
$$h^1(\calm_i)=-\Delta_i+h^0(\calm_i)=-\Delta_i+h^1(\calm_i^{-1}(K+E_{j(i)}))\leq
-\Delta_i+M^{g_{j(i)}}(-\Delta_i).$$

\para\label{UJ4} {\bf Remark.} Assume that we add another  term
$x_{t+1}=x_t+E_{j(t)}$ to the sequence $\{x_i\}_{i=1}^t$ with
similar restriction  $x_{t+1}\in l'-K-L_{\Q,ne}$. Then
$deg_{E_{j(t)}}\calm_t>2g_{j(t)}-2$,  $\Delta_t\geq g_{j(i)}$ and
$M^{g_{j(t)}}(\Delta_t)=h^1(\calm_t)=0$. Therefore, even if one
continues the sequence arbitrarily long inside of
$l'-K-L_{\Q,ne}$, nothing will be changed (e.g. the upper bound
accumulates no more contribution). Sometimes we will just say and
write that $x_t=\infty$, which means that $x_t$ is in the `right'
region $l'-K-L_{\Q,ne}$.

\vspace{2mm}

Next  we  reinterpret (\ref{UJ2}) in terms of path cohomology.
Let $\calP$ be the set of paths with $x_0=0$ and $x_t=\infty$, in
the sense of (\ref{UJ4}). Moreover, consider the weight functions
$\{w_q\}_q$ associated with $(\Gamma,k)$ as in (\ref{dEF1}), and
write $\bH^0(\gamma;\Gamma,k)$ for $\bH^0(\gamma,w)$. Then from
(\ref{UJ2}) and (\ref{PC2}) we get

\para\label{PU2} {\bf Corollary.} {\em For any $\gamma\in
\calP$ one has \ $h^1(\call)\leq eu\, (\bH^0(\gamma;\Gamma,k))$.
Hence}
$$h^1(\call)\leq \min_{\gamma\in\calP}\ eu\, (\bH^0(\gamma;\Gamma,k)).$$

\para\label{PU3} {\bf Remark.} Recall that by  (\ref{PC4}) one has: \
$ \min_{\gamma\in\calP}\ eu\, (\bH^0(\gamma;\Gamma,k))\ \leq
 \ eu\, (\bH^0(\Gamma,k))$.

\para\label{PU4} {\bf Example.} If $\Gamma$ is almost rational
(cf. \ref{AR1}), a consequence of the results of \cite{NOSZ} is that
$$\min_{\gamma\in\calP}\ eu\, (\bH^0(\gamma;\Gamma,k))\, = \,
 eu\, (\bH^0(\Gamma,k)),$$
and, in fact, the minimum $\min_{\gamma\in\calP}$ is realized by
an increasing path. The point is that $\ker{U}\in \bH^0(\Gamma,k)$
admits `representative' lattice points which are  totally ordered
(with respect to $<$) sitting on an increasing path. In fact,
$\bH^0(\Gamma,k)$ is determined in \cite{NOSZ} from the values of
$\chi_k$ along this path.

\para\label{PU5} {\bf Example.} The situation from (\ref{PU4}), in general,  is
not true. I.e.,  one may have
$$ \min_{\gamma\in\calP}\ eu\, (\bH^0(\gamma;\Gamma,k))\, <\,
 eu\, (\bH^0(\Gamma,k)),$$
i.e., the path cohomology may provides a strict better upper bound
for $h^1(\call)$ than the lattice cohomology (cf. \ref{PU3}).
To see this, construct $\Gamma$ with $c=g=0$ as follows.
Let $E$ and $E'$ be two vertices, both with self--intersection
$-14$, and connected by an edge. Attach to both of them two--two
cusps as in (\ref{NV1}--\ref{NV2}). Take $k=K$. Then
$\chi(Z_{min})=\chi(3Z_{min})=-3$ and
$m_K=\chi(2Z_{min})=-4$. By a computation
$$\bH^0(\Gamma,K)=\et^+_{-8}\oplus \et_{-6}(1)^6\oplus
\et_0(1)^2,$$ where the generators in degree zero are 0 and $-K$,
in degree $-4$ is $2Z_{min}$, while in degree $-3$ the cycles
$Z_{min},L,R,L',R',3Z_{min}$. Here, the cycles $L$ and $R$ are
symmetric with respect to the natural symmetry compatible with
$E\leftrightarrow E'$, for both $Z_{min}<L,R<2Z_{min}$, but $L$
and $R$ are not comparable by $<$. Hence, when one travels from
$Z_{min}$ to $2Z_{min}$ by a Laufer type path, then one has to
make a choice (left-right) to pass through $L$ or $R$, but one
doesn't have to touch both of them. The situation is similar with
$L'$ and $R'$ which sit between $2Z_{min}$ and $3Z_{min}$. Hence,
it turns out  that the module for a minimal increasing path (with
end--point at $K$, or at $\infty$) is
$$\bH^0(\gamma_{min},K)=\et^+_{-8}\oplus \et_{-6}(1)^4\oplus \et_0(1)^2,$$
which has $eu$ two less than $\bH^0(\Gamma,K)$.

\para\label{PU6} {\bf Example.} We may ask how sharp is the
topological upper bound  (\ref{PU2}). Although it is not very easy
to provide abundant examples for $h^1(\call)$, for the geometric
genus $p_g:=h^1(\calo_{\tx})$  more examples are available. In
this case, in many graphs the inequality (\ref{PU2}) is optimal,
i.e. the topological upper bound is realized by the $p_g$ of some
analytic structure. Nevertheless, this is not the case all the
time. For the graph $\Gamma$ discussed in (\ref{NV1}), both the
lattice and path cohomologies provide the same upper bound
$p_g\leq 4$ (cf. \ref{PU3}). On the other hand, by a (not simple) line of
arguments, one finds out that there is no analytic structure
supported on this topological type with $p_g=4$ ($p_g=3$ can be
realized by a splice type complete intersection).
The reader may decide if this example is `generic' or
`pathological'.

(Note that $p_g\geq h^1 (Z_{min})=1-\chi(Z_{min})$, hence $p_g\geq 2$
for any analytic structure, while $p_g=3$ for any Gorenstein structure.)

\section{The Seiberg--Witten invariant conjecture.}

\subsection{Line bundles on $\tx$ revisited.}\label{LineBund}

\para\label{LB1} {\bf The bundles $\calo_{\tx}(l')$.}
Start with the data of (\ref{CohCom}) and
assume that {\em  $M$ is rational homology sphere}.
The  `exponential exact sequence' $0\to\Z\to\calo\to
\calo^*\to 0$ \, on $\tx$ induces the exact sequence
\begin{equation*}
0\to H^1(\tx,\calo_{\tx})\to Pic(\tx)\stackrel{c_1}{\to} L'\to 0.
\end{equation*}
For any $l\in L $ one has $c_1(\calo_{\tx}(l))=l$.
Hence  $l\mapsto \calo_{\tx}(l)$ is a group
section of $c_1$ above the subgroup $L$ of $L'$. Since $L'/L$ is
torsion, and $H^1(\calo_{\tx})=\C^{p_g}$ is torsion--free, this can be
extended in unique way to a group section $s:L'\to Pic(\tx)$ of
$c_1$. We write $\calo_{\tx}(l')$ for $s(l')$.

\para\label{LB2}{\bf Relation with coverings.} The next theorem
(\ref{LB3}) illuminates a different aspect of the line bundles
$\calo_{\tx}(l')$. Notice that $\tx\setminus E\approx X\setminus
\{0\}$ has the homotopy type of $M$, hence the abelianization map
$\pi_1(\tx\setminus E)=\pi_1(M)\to H$ defines a regular Galois
covering of $\tx\setminus E$. This has a unique extension $p:Z\to
\tx$ with $Z$ normal and $p$ finite \cite{GR}. The (reduced)
branch locus of $p$ is included in $E$, and the Galois action of
$H$ extends to $Z$ as well. Since $E$ is a normal crossing
divisor, the only singularities that $Z$ might have are cyclic
quotient singularities.

\para\label{LB3}{\bf Theorem.} {\em Consider the finite covering $p:Z\to
\tx$, and set $Q\subset L'$ as in (\ref{2.4}). Then the
$H$-eigenspace decomposition of $p_*\calo_Z$ has the form:
$$p_*\calo_Z=\oplus _{\chi\in \hat{H}}\call_\chi,$$
where $\call_{\theta(h)}=\calo_{\tx}(-l'_e(h))$  for any $h\in H$.
In particular,} $p_*\calo_Z=\oplus_{l'\in Q}\calo_{\tx}(-l').$

\vspace{2mm}

\noindent The proof is based on a similar statement of Koll\'ar
valid for cyclic coverings, see e.g. \cite{Kollar}, \S 9. For
details, see \cite{Line,trieste} or \cite{Okuma}.

\subsection{The conjectured identities.}

\para The next expected
property is a generalization of the conjecture of \cite{[51]},
where only the case of canonical $spin^c$--structure was
considered. The generalization to any $spin^c$--structure appeared
in \cite{Line}, where it was formulated for any $\Q$--Gorenstein
singularity (with rational homology sphere link). The article
\cite{SI} shows that we cannot expect the validity of the
identities in this generality. Nevertheless, we expect that
it is true for a large class of normal surface singularities
(subclass of $\Q$-Gorenstein singularities with rational homology
sphere links). In the next paragraphs we will present two
(equivalent) versions.

In this section we assume that the link $M$
of $(X,0)$ is a rational homology sphere.
We fix a good resolution $\pi:\tx\to X$ with $s:=\#\calj$. Also, we set
$$\bL':=\{ l'\in L':\, e(l')=l'_{ne}(h) \ \mbox{for some $h\in H$}\}=
\bigcup_{h\in H}l'_{ne}(h)+L_e.$$ (For notations, see (\ref{2.4}) and
(\ref{4.4}).) One can verify that
 $L'_e\subset  \bigcup_{l'\in Q} -l'+L_e
\subset \bL'$.

\para\label{5.4s}{\bf Property A.} {\em
 Consider an arbitrary $l'\in
\bL'$ and define a characteristic element by $k:=K-2l'\in Char$.
Then, we say that $(X,0)$ satisfies Property A if}
\begin{equation*}
h^1(\calo_{\tx}(l'))= -\ssw(M,[k])-\frac{k^2+s}{8}.\tag{1}
\end{equation*}

\para\label{5.5}{\bf Remark.} In order to prove the property,
it is enough to verify it for line bundles $\call$ with
$c_1(\call)=l'$ of type $l'=l'_{ne}(h)$ (for some $h\in H$).
Indeed, write $l'$ in the form $l'=l'_1+l$ where
$l'_1=e(l')=l'_{ne}(l'+L)$ and $l\in L_e$. Let $RHS(l')$, resp.\
$RHS(l'_1)$, be the right hand side of (1) for $l'$, resp.\
$l'_1$. Since $[K-2l']=[K-2l'_1]$, the Seiberg-Witten invariants
are the same, hence
$$RHS(l')-RHS(l'_1)=\frac{-(K-2l')^2+(K-2l'_1)^2}{8}=-(l,l')-\chi(l).$$
This combined with (\ref{4.4})(c) shows that
(\ref{5.4s})(1) for $\call$  and $\call\otimes \calo_{\tx}(-l) $ are
equivalent.

In fact, consider {\em any}  set of representatives $\{l'\}_{l'\in
R}$ $(R\subset \bL')$ of the classes $H$, i.e.\ $\{l'+L\}_{l'\in
R}=H$. Then the above argument applied for elements from $R$ shows
that the validity of the property (\ref{5.4s})  follows from the
verification of (1) for line bundles $\call$ with $c_1(\call)\in
R$. The possibility $R=-Q$ is emphasized by (\ref{LB3}) and will
be exploited in the second version of the property.

\para\label{dis}{\bf Universal abelian cover.} Let
$(X_{ab},0)$ be the universal
abelian cover of $(X,0)$ with its natural $H$-action. Namely,
$(X_{ab},0)$ is the unique normal singularity with a finite
projection $(X_{ab},0)\to (X,0)$, regular over $X\setminus 0$
corresponding to the abelianization map  $\pi_1(X \setminus
0)=\pi_1(M)\to H$.  Then the space $Z$ considered in
(\ref{LB2}--\ref{LB3}) is a partial resolution of $(X_{ab},0)$ with
only cyclic quotient singularities. The geometric genus
$p_g(X_{ab},0)$ of $(X_{ab},0)$ can be computed as the dimension
of $H^1(Z,\calo_{Z})$,  but this space has a natural eigenspace
decomposition $\oplus _{\chi\in \hat{H}}H^1 (Z,\calo_{Z})_\chi$ too.
Hence one may consider the invariants
$$p_g(X_{ab},0)_\chi:=\dim_\C H^1(Z,\calo_Z)_\chi \ \  \
\mbox{(for any $\chi\in \hat{H}$)}.$$ Notice that  (\ref{LB3})
reads as
$$p_g(X_{ab},0)_{\theta(h)}=h^1(\calo_{\tx}(-l'_e(h))) \ \ \ \mbox{(for any $h\in H$)}.$$


Since the set $\{-l'_e(h)\}_{h\in H}$ is a set of representatives
for $H$, by (\ref{5.5}) the previous  Property A (\ref{5.4s}) is
{\em equivalent} with the following.

\para\label{conj2}{\bf Property B.} {\em For any $h\in
H$ consider  $k:=K+2l'_e(h)\in Char$. Then for any $h\in H$
\begin{equation*}
p_g(X_a,0)_{\theta(h)}= -\ssw(M,[k])-\frac{k^2+s}{8}. \tag{2}
\end{equation*}

\subsection{Examples.}

 \para\label{6.2}{\bf Example.}} {\bf Property  A}
(hence B  too)  {\bf is true for any rational singularity.}
Indeed,  by (\ref{5.5}), we can assume  that $l'=l'_{ne}(h) $
for some  $h$. Then, by (\ref{4.3}),  $h^1(\calo_{\tx}(l'))=0$.
On the other hand, by \cite{NOSZ}, $-\ssw(M,[k])=(k_r^2+s)/8$,
where
$k_r=K+2\bar{l}'_{ne}(-l'+L)$.  Since
$\bar{l}'_{ne}(-l'+L)=-l'_{ne}(l'+L)=-l'$ one gets
$k_r=k$. Hence the right hand side of
(\ref{5.4s})(1) is also vanishing.

This proof also shows that for $(X,0)$ rational, and for any $h\in
H$, one has
$$p_g(X_{ab},0)_{\theta(h)}=\frac{(K+2\bar{l}'_{ne}(h))^2-(K+2l'_e(h))^2}{8}=
-\chi(\bar{l}'_{ne}(h))+\chi(l'_e(h)).$$ In particular,
$(X_{ab},0)$ is rational if and only if
$\chi(\bar{l}'_{ne}(h))=\chi(l'_e(h))$ for all $h\in H$.
One can find rational graphs whose
universal abelian covers are not rational, a fact which stresses
the differences between the `liftings' $l'_e(h)$ and $\bar{l}'_{ne}(h)$.

\para\label{SQ} {\bf Example. Splice quotients.} The validity of Property A for
rational singularities (cf. \ref{6.2}), the surgery formulas of
\cite{BN} regarding the Seiberg--Witten invariants, and
the result of Okuma from \cite{Okuma2} lead in \cite{BN}  to the
verification of Property A for all splice quotients.
 (The case of trivial line bundle was verified earlier in \cite{NO}.)
Splice quotient singularities were introduced by Neumann and Wahl
(see e.g. \cite{NWuj2}), they include all the rational, minimal
elliptic singularities, and all singularities which admit a good
$\C^*$--action.

Assume now that $(X,0)$ is a splice quotient, and additionally,
its topological type is also almost rational. Set $l'\in \bL'$ and
$k=K-2l'$ as in Property A. Then Property A, (\ref{HF4}) and
(\ref{PU4}) read as
$$h^1(\calo_{\tx}(l'))=eu(\bH^0(\Gamma,k))=\,
\min_{\gamma\in\calP}\, eu(\bH^0(\gamma;\Gamma,k)),$$ which
(by \ref{PU2}) is a topological upper bound for  $h^1(\call)$,
where $\call$ is an {\em any bundle} with
$c_1(\call)=l'$.

In particular, if $\Gamma$ is almost rational, and the topological
type admits a splice quotient analytic structure, then the
geometric genus of the splice quotient analytic structure (which
satisfies Property A) is an upper bound for the geometric genera
of all the possible analytic structures.


\para\label{EX2} {\bf Example.} One can find
even hypersurface singularities when Property A is not true
for $p_g$ (i.e. for $l'=0$). Such examples are provided in \cite{SI} by
super--isolated singularities. In the examples of \cite{SI}(4.1),
$p_g$ is strict higher then the expected  value
$-\ssw(M,[K)]-(K^2+s)/8$.
Now, using our  previous discussions, this phenomenon
can be explained as follows. 

In general, in the light of Conjecture (\ref{HF5}),
Property A/B is equivalent to
\begin{equation*}
p_g=eu(\bH^*(\Gamma,k))=-\min\chi_{K}+\textstyle{\sum_q}(-1)^q
\rank_{\Z}\bH_{red}^q(\Gamma,K).\tag{1}
\end{equation*}
On the other hand, the inequalities from subsection (\ref{PCUB})
read as
\begin{equation*}
p_g\leq \min_{\gamma\in\calP}\ eu\,
(\bH^0(\gamma;\Gamma,K))\,\leq  \, eu(\bH^0(\Gamma,K))=
-\min\chi_{K}+ \rank_{\Z}\bH_{red}^0(\Gamma,K).\tag{2}
\end{equation*}
Assume that three things are happening simultaneously: (a) in (2)
the second inequality is equality, (b) for some analytic structure
the first inequality in (2) is sharp (hence
$p_g=eu(\bH^0(\Gamma,K))$), and (c) $\bH^q_{red}\not= 0$
for $q\geq 2$, creating the situation  $eu(\bH^0)>eu(\bH^*)$. Then
Property A fails, and in fact $p_g>eu(\bH^*(\Gamma,k))$ for that
analytic structure.

This is the case for all the examples of (4.1)\cite{SI}.

Let us analyse a little bit more the case $C_4$ of \cite{SI}. The
corresponding graph $\Gamma$ is

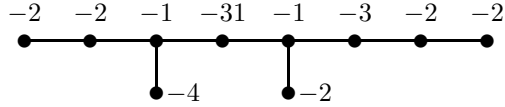
\begin{picture}(300,45)(20,0)
\put(125,25){\circle*{5}} \put(150,25){\circle*{5}}
\put(175,25){\circle*{5}} \put(200,25){\circle*{5}}
\put(225,25){\circle*{5}} \put(150,5){\circle*{5}}
\put(200,5){\circle*{5}} \put(100,25){\line(1,0){175}}
\put(150,25){\line(0,-1){20}} \put(200,25){\line(0,-1){20}}
\put(125,35){\makebox(0,0){$-2$}}
\put(150,35){\makebox(0,0){$-1$}}
\put(175,35){\makebox(0,0){$-31$}}
\put(200,35){\makebox(0,0){$-1$}}
\put(225,35){\makebox(0,0){$-3$}} \put(160,5){\makebox(0,0){$-4$}}
\put(210,5){\makebox(0,0){$-2$}} \put(100,25){\circle*{5}}
\put(250,25){\circle*{5}} \put(275,25){\circle*{5}}
\put(100,35){\makebox(0,0){$-2$}}
\put(250,35){\makebox(0,0){$-2$}}
\put(275,35){\makebox(0,0){$-2$}}
\end{picture}

In this case $\bH^0(\Gamma,K)=\et_{-10}^+\oplus \et_{-10}(3)\oplus
\et_0(1)^2$, hence $eu(\bH^0)=10$, but $eu(\bH^*)=8$.
(Strictly speaking, the author verified that $-\ssw(M,[K])-(K^2+s)/8=8$,
cf. (\ref{HF5}).) 
 Hence the
topological bound given by (2) is $p_g\leq 10$.   This topological
type admits two, very natural, but rather different analytic
structures. The first is the super--isolated hypersurface
singularity mentioned above: it has $p_g=10$ \cite{SI}. On the
other hand, there is a splice quotient singularity which satisfies
Property A, hence with $p_g=8$ \cite{NO}. This is the
$\Z_5$--factor of the complete intersection
$\{z_1^3+z_2^4+z_3^5z_4=z_3^7+z_4^2+z_1^4z_2=0\}\subset (\C^4,0)$
by the diagonal action $(\alpha^2,\alpha^4,\alpha,\alpha)$
($\alpha^5=1$).

Therefore, in general, the geometric genus of those analytic
structures which satisfy Property A is not `extremal' (in contrast
with the almost rational case (\ref{SQ})). In \cite{BN},  Property
A is reformulated completely in terms of the analytic structure
(independently of any Seiberg--Witten type theory).  \cite{Okuma2}
suggests that in the heart of the its validity there is a
cohomological vanishing result.

}

\end{document}